\documentclass[11pt,a4paper,twoside]{article}
\usepackage[top=2.3cm,bottom=2.3cm,left=2.4cm,right=2.4cm,headheight=14pt]{geometry}
\usepackage[T1]{fontenc}
\usepackage{amsmath,amsthm,mathtools}
\usepackage{mathptmx}
\usepackage{amssymb}
\usepackage{courier}
\usepackage{microtype,booktabs,threeparttable,graphicx,xcolor,enumitem,tabularx,caption,makecell}
\usepackage{algorithm}
\usepackage[noend]{algpseudocode}
\usepackage{setspace}
\usepackage[round,authoryear,sort]{natbib}
\usepackage[hidelinks]{hyperref}
\hypersetup{pdftitle={Forecast Projection and Decision Regret: Theory and Evidence from Battery Arbitrage},pdfauthor={Debangana Mukherjee and Akash Ashirbad Panda}}
\usepackage{orcidlink}
\usepackage{fancyhdr}
\setlist{topsep=2pt,itemsep=1pt,parsep=0pt}
\newtheorem{theorem}{Theorem}
\newtheorem{proposition}{Proposition}
\newtheorem{lemma}{Lemma}
\newtheorem{corollary}{Corollary}
\theoremstyle{definition}\newtheorem{definition}{Definition}\newtheorem{assumption}{Assumption}
\theoremstyle{remark}\newtheorem{remark}{Remark}
\newcommand{\1}{\mathbf 1}
\newcommand{\R}{\mathcal R}
\newcommand{\F}{\mathcal F}
\newcommand{\E}{\mathbb E}
\DeclareMathOperator*{\argmax}{arg\,max}
\DeclareMathOperator*{\argmin}{arg\,min}
\newcommand{\IEXnDays}{1{,}522}
\newcommand{\IEXoracle}{65.82}

\newcommand{\IEXvOneVal}{6.55}
\newcommand{\IEXvThreeVal}{4.86}

\newcommand{\IEXvOneLoss}{1.69}
\newcommand{\etaKa}{324}
\newcommand{\etaKd}{2{,}496}

\newcommand{\IEXmlpVsEw}{$+0.47$}

\newcommand{\DEbestVsEw}{$+167$}
\newcommand{\DEbestVsEwCI}{$-$9, 346}
\newcommand{\DEpmin}{$-222$}
\newcommand{\DEpmax}{210}
\newcommand{\DEoracle}{9{,}541}
\newcommand{\DEoracleH}{9{,}256}
\newcommand{\DEsel}{V1 $(\lambda,\text{lr})=(0.1,3\times10^{-3})$; V3 $(\lambda,\text{lr})=(0.1,3\times10^{-3})$; V4 $(\lambda,\text{lr})=(1,3\times10^{-3})$}
\newcommand{\DEshiftVal}{3.1}
\newcommand{\DEvOneH}{2{,}689}
\newcommand{\DEvThreeH}{2{,}702}
\newcommand{\DEewH}{2{,}257}
\newcommand{\DEmlpH}{2{,}425}
\newcommand{\RobustSig}{8}
\newcommand{\ConvLvlA}{2{,}449}
\newcommand{\ConvRegA}{6.60}

\newcommand{\ConvLvlD}{427}
\newcommand{\ConvRegD}{4.58}
\newcommand{\ConvRegDcor}{4.57}

\newcommand{\FreshN}{53}
\newcommand{\FreshOracle}{86.9}
\newcommand{\FreshEW}{2.90}
\newcommand{\FreshReg}{2.67}

\newcommand{\FreshVfourDiff}{$+0.42$}
\newcommand{\FreshVfourCI}{$-0.45$, $+1.35$}
\newcommand{\FreshProjCost}{$+0.15$}
\newcommand{\FreshProjCostCI}{$0.00$, $+0.33$}
\newcommand{\FreshVtwoAbove}{55}
\newcommand{\FreshVtwoLevel}{+230}
\newcommand{\FreshVoneLevel}{-93}

\newcommand{\BudgetEWval}{4.08}
\newcommand{\BudgetEWhold}{2.43}
\newcommand{\BudgetMinGap}{0.67}

\title{Forecast Projection and Decision Regret:\\Theory and Evidence from Battery Arbitrage}
\author{Debangana Mukherjee\,\orcidlink{0000-0002-1516-840X}\thanks{School of Interwoven Arts and Sciences, Krea University, Sri City, India (\href{mailto:debangana.mukherjee@krea.edu.in}{\texttt{debangana.mukherjee@krea.edu.in}}).}
\and Akash Ashirbad Panda\,\orcidlink{0000-0003-3496-9111}\thanks{Corresponding author. Department of Mathematics, Indian Institute of Technology Indore, Simrol, India (\href{mailto:akashpanda@iiti.ac.in}{\texttt{akashpanda@iiti.ac.in}}).}}
\date{}

\begin{document}
\maketitle
\vspace{-1.5em}
\begin{abstract}\noindent
Decision-focused learning trains forecasts for the decisions they lead to, but forecasts are often
projected onto known bounds before use. We study the resulting mismatch in battery arbitrage. When a
battery starts and ends the day at the same charge level, a common shift in prices affects its schedule
only through round-trip losses. As a result, decision-focused losses identify the forecast level only
weakly. For the SPO+ and perturbed Fenchel--Young losses, we characterize the preferred level by a
throughput-matching condition. Projecting forecasts onto a price cap can then remove profitable
spreads and increase regret, even as prediction error falls. We give bounds on the change in regret, a
dual-threshold condition under which a raw-optimal schedule remains optimal after clipping, and a general
version of the mechanism for other problems with a conservation constraint. On 1{,}158 holdout days of the
Indian day-ahead market, an equal-weight average of three statistical forecasts outperforms the five
decision-focused variants tested. On \FreshN\ later days, with 62\% of price blocks at the cap, no variant
outperforms the regression or the equal-weight average, although most differences from the average are not
significant. The results suggest that forecast post-processing should be treated as part of the
predict-then-optimize pipeline.



\end{abstract}

\section{Introduction}\label{sec:intro}
Decision-focused learning (DFL) trains a forecast to improve the decision obtained from an optimization model
rather than to minimize prediction error alone \citep{elmachtoub2022predict,mandi2024dfl}, including in
energy storage \citep{sang2022decision,smets2025forecasting,smets2025value,smets2026proxies}. A forecast, however, need not be used in its raw form: it may first be projected onto known price, demand or capacity bounds; more generally, forecasts are often transformed by clipping, truncation, censoring, rounding or normalization before they are used. We ask whether a forecast direction that is nearly irrelevant to the training loss can become
important to the decision after this projection.

In electricity markets such bounds are often regulatory. The ordinary day-ahead market of the Indian Energy
Exchange had a ceiling of Rs\,12/kWh from April 2022, lowered to Rs\,10/kWh from 4~April 2023
\citep{cerc2022sm4,cerc2022sm5,cerc2023sm4}. Cutting forecast prices above the cap down to the cap before
scheduling is a natural choice, because actual prices cannot be higher and the change can only reduce the
error in each block; the cap does not, however, force the operator to do it.

Consider a price-taking battery with feasible schedules $z\in\F$, net discharge $q(z)$, operating cost
$C(z)$, and realized price vector $p\in[\underline p,\bar p]^H$. For block length $\Delta$, its profit and
perfect-information value are
\begin{equation}\label{eq:intro-profit}
J(p,z)=\Delta p^\top q(z)-C(z),\qquad \Pi^\star(p)=\max_{z\in\F}J(p,z).
\end{equation}
Given a raw forecast $u$, the operator uses the box projection $g(u)=\Pi_{\mathcal B}(u)$ onto
$\mathcal B=[\underline p,\bar p]^H$ and selects $z_g(u)\in\argmax_{z\in\F}J(g(u),z)$. The relevant loss in
operation is the deployed regret
\begin{equation}\label{eq:intro-regret}
\R^g(p,u)=\Pi^\star(p)-J\bigl(p,z_g(u)\bigr),
\end{equation}
which need not follow the prediction error of $g(u)$, nor the decision-focused loss evaluated at $u$.

The key mechanism is a conservation identity. If the battery ends the day at its initial charge, has
round-trip efficiency $\eta$, and charges $T(z)=\Delta\sum_t c_t$ MWh, then
\begin{equation}\label{eq:intro-shift}
\Delta\1^\top q(z)=-(1-\eta)T(z),\qquad J(p+k\1,z)=J(p,z)-k(1-\eta)T(z).
\end{equation}
Common price shifts therefore leave the raw optimization problem unchanged when $\eta=1$; when $\eta<1$, they
alter the profit of each feasible schedule only through its charged throughput. For the SPO+ loss, the
expected loss along $u+k\1$ is convex, and its minimizing shifts satisfy the \emph{throughput-matching}
condition
\begin{equation}\label{eq:intro-matching}
\E\,T^-(2u-p+2k\1)\ \le\ \E\,T(z^\star)\ \le\ \E\,T^+(2u-p+2k\1),
\end{equation}
where the expectation is over forecasts and prices, $T^-$ and $T^+$ are the minimum and maximum charging
throughputs of optimal helper schedules, and $z^\star$ is optimal at the realized prices. Under an explicit
excess-throughput condition (Assumption~\ref{ass:over}), this shift is positive; the same mechanism has a
perturbed Fenchel--Young counterpart (Theorem~\ref{thm:bias}). Squared-error regularization constrains the
population-optimal shift (Corollary~\ref{cor:reg}), whereas incomplete training can create a much larger
out-of-sample level error: about Rs\,2{,}450/MWh with 150 training steps per refit in our data, falling to
about 430 with 5{,}000 steps (Section~\ref{sec:results}).

Projection breaks the near-invariance in \eqref{eq:intro-shift}. Although $|g(u)_t-p_t|\le|u_t-p_t|$ for
every block, clipping can remove the price spreads that determine profitable trades. For example, with
realized prices $(5,10)$, forecast $(7,14)$, cap $10$, no energy loss and a total cost of buying and selling
of $4$, clipping reduces the mean squared error from $10$ to $2$ but changes realized trading profit from
$1$ to $0$. More generally, with no energy loss, a raw forecast can have zero unregularized SPO+ loss while
its projected schedule loses the entire available arbitrage profit (Theorem~\ref{thm:eta1}).

\subsection{Related literature}\label{sec:lit}
\citet{elmachtoub2022predict} proposed the
SPO+ loss for problems with linear objectives. \citet{berthet2020learning} and \citet{blondel2020learning}
proposed perturbed optimizers and Fenchel--Young losses. \citet{mandi2024dfl} review the area. Other work
asks when it is better to train for accuracy first and optimize later. \citet{hu2022fast} show that this
simple approach can work very well when the regression is accurate, and \citet{elmachtoub2023estimate}
compare the two approaches. \citet{bertsimas2020predictive} study how to use data directly for decisions.
\citet{lee2026irrelevant} treat prediction directions that do not affect decisions with a geometric
approach that projects the prediction error onto a tangent space. In our problem the forecast level is such a direction when $\eta=1$, and
nearly so when $\eta<1$; anchoring removes it from training in a similar spirit. Our work adds three
things for this setting. First, when $\eta<1$ the direction is only \emph{nearly} irrelevant, and the loss
has a definite preferred value along it, which we characterize (throughput matching). Second, this preferred value moves away from the accurate level; for the pure loss the shift grows like $1/(1-\eta)$ if the matching cost increase stays bounded away from zero, and with a squared-error term it stays bounded. Third, a separate step after the forecast, projection to
the price cap, makes this direction matter for decisions.

\citet{sang2022decision} train price forecasts for battery trading with a
decision-focused loss. \citet{smets2025forecasting,smets2025value,smets2026proxies} study value-oriented
and proxy-based DFL for batteries. \citet{yi2025perturbed} use perturbed optimizers for battery bidding,
and \citet{toubeau2021datadriven} schedule batteries in energy and reserve markets.
\citet{alkhulaifi2026dfl} report strong results for DFL in energy storage with automated feature
engineering, in a different setting from ours; our contrasting results suggest that the value of DFL
depends on the market and on how the forecast is used. \citet{matsumoto2026guard} show that forecast
noise leads to excessive battery trading and that conservative safety margins (guard bands) reduce this
optimizer's curse \citep{smith2006optimizer}. This is a closely related economic explanation of the overtrading
in Assumption~\ref{ass:over}. Our contribution is different: we characterize exactly where the SPO+ and
Fenchel--Young losses place the forecast level (throughput matching), show that this works like an
implicit safety margin (Remark~\ref{rem:guard}), and study how it interacts with projection to a price
cap. As far as we know, none of these studies looks at price caps, at
how well the forecast level is learned, or at comparisons with forecast averages.

\citet{weron2014electricity} and
\citet{lago2021forecasting} review price forecasting and propose the LEAR and DNN benchmarks.
\citet{maciejowska2026economic} argue that forecasts should be judged by their economic value, not only
by accuracy. Work on combining forecasts \citep{bates1969combination,timmermann2006forecast} shows that
equal weights are hard to beat \citep{smith2009simple}. We find the same for decisions: weights chosen from past profit do not outperform equal weights.

\citet{sioshansi2009pjm}, \citet{antweiler2021storage} and
\citet{mercier2023storage} measure the value of battery trading and how it depends on forecasts and
market rules. We add evidence on how a price cap changes the value of different forecasting
\emph{methods}.

\subsection{Contributions}\label{sec:contrib}
Neither the convexity of SPO+ nor the energy balance of a battery is new. What is new is how they interact
with a step that changes the forecast after training and before dispatch.

On the theory side, Theorem~\ref{thm:bias} gives the price level that minimizes the expected SPO+ loss through
the matching condition \eqref{eq:intro-matching}. It shows when this level lies above the accurate level and
how it depends on efficiency, and Corollary~\ref{cor:reg} shows how a squared-error term limits it.
Theorem~\ref{thm:general} states the same mechanism for any decision problem with
$\Delta\,r^\top q(z)=-a\,T(z)+b$ and bounded $T$, for example a flexible load that must buy a fixed amount of
energy. \citet{lee2026irrelevant} study prediction directions that do not affect decisions; we study what
happens when the forecast is changed after training.

We also study regret after projection. A raw forecast can have zero SPO+ loss and still have positive regret
once it is clipped (Theorem~\ref{thm:eta1}), so a regret bound for the raw loss does not carry over to the
deployed forecast. Proposition~\ref{prop:margin} gives the exact condition, in the form of a normal cone, under
which any post-processing step keeps a raw-optimal decision optimal, and bounds the change in regret when it
does not. Theorem~\ref{thm:safe} gives a simpler sufficient condition for batteries that can be read off the
dual solution of the raw problem, and Corollary~\ref{cor:select} deals with ties between equally good
schedules. These results add to storage guard-band methods \citep{matsumoto2026guard} rather than replace them.

We then propose a simple fix. Taking the forecast level from a separate model, which we call anchoring, stops
training from moving the level (Corollary~\ref{cor:anchor}). In exploratory comparisons
(Section~\ref{sec:rule}), anchoring decision-focused forecasts every day worked better than anchoring only when
the certificate failed, statistical forecasts did best with their own level, and keeping the raw schedule
whenever it remained optimal made the deployed regret equal to the raw regret.

Finally, we test the mechanism on 1{,}158 holdout days of the Indian day-ahead market and on \FreshN\ later
days in which most prices were at the cap, and compare with a German market in which the limits do not bind
(Section~\ref{app:de}). An equal-weight average of three statistical forecasts has the lowest regret on the
long holdout, 0.47 [0.34, 0.59] lakh per MW-year less than the best DFL variant, or about 0.7\% of the
perfect-foresight profit. The fresh test gives no sign that any of the tested DFL variants does better. These
findings are about the DFL implementations we tested, not about DFL methods in general.

The rest of the paper is organized as follows. Section~\ref{sec:model} describes the battery dispatch problem, projection and regret. Section~\ref{sec:theory} gives the theory.
Section~\ref{sec:design} describes the data and test design, Section~\ref{sec:results} gives the results,
and Section~\ref{sec:impl} discusses what they mean and their limits. The electronic companion contains the proofs (Section~\ref{app:proofs}), implementation details
(Section~\ref{app:impl}), the German test (Section~\ref{app:de}), further results (Section~\ref{app:extra}) and the details of
the operating rules (Section~\ref{app:rules}).

\section{Battery dispatch, projection and regret}\label{sec:model}
A delivery day has $H$ blocks, each $\Delta$ hours long (for the IEX,
$H=96$ and $\Delta=0.25$). Prices $p\in\mathbb R^H$ must lie between a floor and a cap set by the
regulator, $p\in[\underline p,\bar p]^H$. The battery is small, so it does not affect prices. It has
energy capacity $B$, power limits $\bar c$ (charging) and $\bar d$ (discharging), and charging and
discharging efficiencies $\eta_c$ and $\eta_d$. The round-trip efficiency is $\eta=\eta_c\eta_d$. A
schedule $z=(c,d,x)$ must satisfy
\[
x_t=x_{t-1}+\Delta\big(\eta_cc_t-d_t/\eta_d\big),\quad 0\le c_t\le\bar c,\quad 0\le d_t\le\bar d,\quad
c_td_t=0,\quad 0\le x_t\le B,\quad x_0=x_H=\bar x .
\]
Let $\F$ be the set of such schedules. The condition $c_td_t=0$ says that the battery cannot charge and
discharge in the same block. The battery starts and ends the day at the same level $\bar x$.
Buying and selling each pay an exchange fee and a battery wear cost $\gamma$ per MWh. The total cost per
MWh is $\omega_b=\phi_b+\gamma$ for buying and $\omega_s=\phi_s+\gamma$ for selling. At prices $\pi$ the
schedule earns
\[
J(\pi,z)=\Delta\textstyle\sum_t\big[(\pi_t-\omega_s)d_t-(\pi_t+\omega_b)c_t\big]=\Delta\pi^\top q(z)-C(z),
\]
where $q=d-c$ is the net energy sold in each block and $C(z)=\Delta\sum_t(\omega_sd_t+\omega_bc_t)$ is the
cost. Let $\Pi^\star(\pi)=\max_{z\in\F}J(\pi,z)$ be the best profit at prices $\pi$, and let
$T(z)=\Delta\sum_tc_t$ be the energy charged during the day.

Without the condition $c_td_t=0$, the problem is a linear
program (LP). The LP can be wrong when prices are very negative: by charging and discharging in the same
block, the battery buys energy and loses part of it in conversion, without changing its stored energy.
Let
\begin{equation}\label{eq:theta}
\theta=-\frac{\omega_b+\eta\,\omega_s}{1-\eta}.
\end{equation}

\begin{lemma}[Exact LP]\label{lem:lp}
If $\pi_t>\theta$ in every block, every optimal solution of the LP has $c_td_t=0$ in every block, so it is
also optimal for the problem with $c_td_t=0$. If $\pi_t<\theta$ in some block, charging and discharging in
that block can increase LP profit. If $\eta=1$ and $\omega_b+\omega_s>0$, the LP is exact at every price.
\end{lemma}

We therefore solve the LP and, if a solution charges and discharges in the same block, solve the problem again with binary variables that rule this out. The LP relaxation is therefore exact only under the price condition of Lemma~\ref{lem:lp}; when helper or training prices can fall below $\theta$, the exact re-solve is needed. With our costs, $\theta$ is about $-6{,}780$ Rs/MWh in India and $-77$ EUR/MWh in Germany. Section~\ref{app:impl} reports how often the exact re-solve was
needed. All results below hold for $\F$ with the condition $c_td_t=0$: the proofs use only that profit is affine in prices and linear in the schedule, that $\F$ does not depend on prices, and that the convex hull of $\F$ has finitely many extreme points. We write $\mathcal V=\operatorname{ext}(\operatorname{conv}\F)$ for this finite set. Because $\F$ is compact, $\mathcal V\subseteq\F$, and a profit that is linear in the schedule attains its maximum over $\F$ at a point of $\mathcal V$; every optimization over $\F$ below can therefore be restricted to $\mathcal V$. In particular, $\Pi^\star(\pi)=\max_{z\in\mathcal V}J(\pi,z)$ is still convex, as a maximum of finitely many affine functions of prices.

On day $D-1$ the operator makes a forecast $\hat p$ using the information
$\mathcal I_D$ available before bidding. Any forecast price outside the allowed range is moved to the
nearest limit, in each block:
\begin{equation}\label{eq:proj}
\Pi_{\mathcal B}(\hat p)=\min\big\{\max\{\hat p,\underline p\},\bar p\big\}.
\end{equation}
We call this step projection. We distinguish two kinds of regret.
\begin{itemize}[leftmargin=1.2em]
\item \emph{Raw regret.} If the operator commits to a schedule $\hat z^{\,\mathrm{raw}}\in\argmax_{z\in\F}
J(\hat p,z)$ that is optimal for the forecast itself, its raw regret is
\begin{equation}\label{eq:rawregret}
\R^{\mathrm{raw}}(p,\hat p)=\Pi^\star(p)-J(p,\hat z^{\,\mathrm{raw}})\ge0.
\end{equation}

\item \emph{Regret after projection.} If the operator first projects the forecast and then commits to a
schedule $\hat z\in\argmax_{z\in\F}J(\Pi_{\mathcal B}(\hat p),z)$, the regret is
\begin{equation}\label{eq:clipregret}
\R^{\mathrm{clip}}(p,\hat p)=\R^{\mathrm{raw}}\big(p,\Pi_{\mathcal B}(\hat p)\big)=\Pi^\star(p)-J(p,\hat z).
\end{equation}

\end{itemize}
Both are paid at the actual prices and measure the profit lost compared with knowing the prices in
advance. Unless stated otherwise, ``regret'' means regret after projection, $\R^{\mathrm{clip}}$, which is what
we report. When several schedules are equally good, a fixed rule picks one. Unless a specific rule is stated
(Corollary~\ref{cor:select} and Section~\ref{sec:rule}), the results below hold for any such rule.

An operator has two reasons to project a forecast. First, realized prices always lie in
$\mathcal B=[\underline p,\bar p]^H$, so projection moves every forecast price weakly closer to the realized one,
$|\Pi_{\mathcal B}(\hat p)_t-p_t|\le|\hat p_t-p_t|$ for every block $t$, and any error measure that increases
with the absolute errors, such as the mean squared or mean absolute error, can only fall. Second, a forecast
above the cap predicts a price that cannot occur. Projection is therefore a natural default, and it looks
harmless. It is still a modeling choice, not a requirement: the cap restricts realized prices, not the prices
an operator feeds into its dispatch model. Our results concern decision-focused forecasts that are projected
before dispatch, not every use of DFL in a capped market, and Section~\ref{sec:results} also reports results
without projection.

\textbf{A two-period example.} Projection can reduce forecast error and still make the decision worse.
Take a battery with no energy loss that can buy 1~MWh in period 1 and sell it in period 2. The total
buying and selling cost is 4 per MWh. Actual prices are $p=(5,10)$, the cap is 10, and the forecast is
$\hat p=(7,14)$ (all in the same price unit).
\begin{itemize}[leftmargin=1.2em]
\item \emph{Before projection}, the forecast spread is $14-7=7>4$, so the battery trades. It earns
$10-5-4=1$.
\item \emph{After projection}, the forecast is $(7,10)$ and the spread is $3<4$, so the battery does not
trade. It earns 0.
\end{itemize}
Projection lowers the mean squared error of the forecast from 10 to 2, but lowers profit from 1 to 0. The
forecast is too high on average, like the DFL forecasts in Section~\ref{sec:results}, and projection removes
the information that period 2 is much more expensive than period 1. This is the main mechanism of the
paper, without SPO+ or a 96-block schedule. It also shows that a smaller forecast error does not always
give a better decision.

\begin{lemma}[Energy balance and price level]\label{lem:cons}
For every $z\in\F$, $\Delta\1^\top q(z)=-(1-\eta)\,T(z)$. Therefore, for all $\pi$ and $k\in\mathbb R$,
\[
J(\pi+k\1,z)=J(\pi,z)-k(1-\eta)T(z)=J\big(\pi,z;\ \omega_b+k(1-\eta),\ \omega_s\big),
\]
and for $\alpha>0$,
\begin{equation}\label{eq:scaling}
\argmax_zJ(\alpha\pi+k\1,z;\omega_b,\omega_s)=\argmax_zJ\big(\pi,z;(\omega_b+k(1-\eta))/\alpha,\ \omega_s/\alpha\big).
\end{equation}

\end{lemma}

Lemma~\ref{lem:cons} says two things. First, adding the same amount $k$ to all prices has exactly the
same effect as raising the cost of buying by $k(1-\eta)$. When $\eta=1$ it has no effect at all. Second,
scaling and shifting a forecast is the same as changing the buying and selling costs in the LP. So
safety margins \citep{matsumoto2026guard} and level changes are the same kind of correction.

\begin{proposition}[Which forecast errors matter]\label{prop:bound}
Let $v=\Pi_{\mathcal B}(\hat p)$ and let $e=v-p$ be the error of the projected forecast, with
daily mean $k$ and remaining part $e^\circ=e-k\1$. Let $\hat z$ be optimal at $v$. Then
\[
0\le\R^{\mathrm{clip}}(p,\hat p)\le\Delta\sum_t|e_t|\,|q^\star_t-\hat q_t|
\quad\text{and}\quad
\R^{\mathrm{clip}}(p,\hat p)\le|k|(1-\eta)\,|T(z^\star)-T(\hat z)|+\Delta\sum_t|e^\circ_t|\,|q^\star_t-\hat q_t|.
\]
\end{proposition}

Proposition~\ref{prop:bound} shows that an error in the forecast level costs money only through the
round-trip loss $1-\eta$ and the difference in charged energy between the two schedules. An error in the
shape of the daily curve costs money directly. The forecast level therefore matters little for the decision, as
long as nothing changes it later. Section~\ref{sec:theory} shows that it starts to matter when the
forecast is trained on a decision loss \emph{and} then projected.

\begin{proposition}[Regret that can and cannot be reduced]\label{prop:decomp}
Let $\nu=\E[p\mid\mathcal I_D]$. For any $\mathcal I_D$-measurable forecast with committed schedule
$\hat z$,

\begin{equation}\label{eq:decomp}
\E\,\R^{\mathrm{clip}}(p,\hat p)=\E\big[\Pi^\star(p)-\Pi^\star(\nu)\big]+\E\big[\Pi^\star(\nu)-J(\nu,\hat z)\big],
\end{equation}
and both terms are nonnegative.
\end{proposition}

The first term is the value of information that no forecast can recover. The second term comes from the
forecasting method, and it is the only part that a better method can reduce. All our comparisons of
methods are differences in this second term.

\section{Decision-focused learning and forecast post-processing}\label{sec:theory}
This section studies how decision-focused losses treat the forecast level and how projection changes the
resulting decisions. Section~\ref{sec:losses} defines the two losses, and Section~\ref{sec:general} states the
mechanism for a general decision problem. Sections~\ref{sec:noloss} to~\ref{sec:training} then develop the
battery case: the lossless battery, the level shift under energy loss, conditions under which projection
preserves a decision, and training on projected forecasts.

\subsection{Decision-focused losses}\label{sec:losses}
\begin{definition}[Decision-focused losses]\label{def:losses}
Fix a perfect-foresight schedule $z^\star\in\argmax_{z\in\F}J(p,z)$. The \emph{SPO+ loss}
\citep{elmachtoub2022predict}, written for our profit-maximization problem, is
\[
\ell(\hat p,p)=\Pi^\star(2\hat p-p)-2J(\hat p,z^\star)+\Pi^\star(p),
\]
with subgradient $2\Delta\big(q(\tilde z)-q(z^\star)\big)$, where $\tilde z$ is the best schedule at the
\emph{comparison price} $2\hat p-p$. The \emph{perturbed Fenchel--Young loss}
\citep{berthet2020learning,blondel2020learning} with noise scale $\sigma>0$ is
\begin{equation}\label{eq:fy}
L_\sigma(\hat p,p)=\E_Z\,\Pi^\star(\hat p+\sigma Z)-J(\hat p,z^\star),
\end{equation}
where $Z$ is standard normal. Its
gradient is $\Delta\big(\E_Zq(\hat z_\sigma)-q(z^\star)\big)$.
\end{definition}

Both losses are convex in $\hat p$. SPO+ satisfies $\ell(p,p)=0$ and $\ell(\hat p,p)\ge\R^{\mathrm{raw}}
(p,\hat p)$: it bounds the regret of the schedule optimized at the \emph{same} forecast that enters the
loss (Section~\ref{app:proofs}). Applied to the projected forecast, this gives
$\ell(\Pi_{\mathcal B}(\hat p),p)\ge\R^{\mathrm{clip}}(p,\hat p)$. In general, $\ell(\hat p,p)$ does \emph{not} bound the
regret after projection, $\R^{\mathrm{clip}}(p,\hat p)$ (Remark~\ref{rem:counter}). In both cases, the gradient compares the schedule of a \emph{helper} dispatch
problem with the perfect-foresight schedule. This comparison drives the results below.

\begin{remark}[The raw SPO+ loss can be zero while regret after projection is positive]\label{rem:counter}
Take the two-period example of Section~\ref{sec:model}: $p=(5,10)$, $\hat p=(7,14)$, cap 10, total trading
cost 4, no energy loss. The perfect-foresight schedule buys in period 1 and sells in period 2, so
$\Pi^\star(p)=1$. At the comparison price $2\hat p-p=(9,18)$ the best profit is $18-9-4=5$, and
$J(\hat p,z^\star)=14-7-4=3$. So
\[
\ell(\hat p,p)=5-2\cdot3+1=0,\qquad \R^{\mathrm{raw}}(p,\hat p)=0,\qquad \R^{\mathrm{clip}}(p,\hat p)=1.
\]
The raw loss is zero, but the deployed (projected) forecast loses all the profit. The bound for the
projected forecast does hold: at $\Pi_{\mathcal B}(\hat p)=(7,10)$ the comparison price is $(9,10)$, with best
profit 0, and $J((7,10),z^\star)=-1$, so $\ell(\Pi_{\mathcal B}(\hat p),p)=0+2+1=3\ge1$.
Theorem~\ref{thm:eta1}(ii) below gives the same failure for a whole day of blocks.
\end{remark}

\subsection{The general mechanism}\label{sec:general}
We first state the mechanism in general: a forecast feeds a decision problem and is post-processed before
use. The rest of this section develops the battery case, which is its main instance and for which the
structure of storage gives sharper results. Let $\mathcal Z$ be a compact feasible set whose convex hull has finitely many extreme points, $\mathcal V_{\mathcal Z}=\operatorname{ext}(\operatorname{conv}\mathcal Z)$ (for example, a polytope, a finite union of polytopes or a bounded mixed-integer set), and let a forecast $u\in\mathbb R^H$ give the profit
$J(u,z)=\Delta\,u^\top q(z)-C(z)$, with $q$ and $C$ linear in $z$. Write $\Pi^\star(u)=\max_{z\in
\mathcal Z}J(u,z)$, and define the SPO+ and perturbed Fenchel--Young losses as in
Definition~\ref{def:losses}. The forecast is deployed after a post-processing map $g:\mathbb R^H\to\mathbb R^H$, such as clipping or projection onto bounds, truncation, censoring, rounding or normalization, and its regret is $\R^g(p,u)=\Pi^\star(p)-J(p,\hat z)$ with
$\hat z$ optimal at $g(u)$.

\begin{assumption}[Near conservation]\label{ass:cons}
There are a direction $r\in\mathbb R^H$, $r\ne0$, a constant $a\ge0$, a constant $b$ and a linear function $T$
with values in $[0,\bar T]$ on $\mathcal Z$, where $\bar T=\max_{\mathcal Z}T$, such that
\begin{equation}\label{eq:cons}
\Delta\,r^\top q(z)=-a\,T(z)+b\qquad\text{for every }z\in\mathcal Z.
\end{equation}

\end{assumption}

\begin{theorem}[Weakly identified directions and post-processing]\label{thm:general}
Under Assumption~\ref{ass:cons}:
\begin{enumerate}[label=(\roman*),leftmargin=2em]
\item $|\ell(u+kv,p)-\ell(u,p)|\le2|k|a\bar T$ for all $u,p,k$. If $a=0$, the SPO+ loss and
the raw decision do not depend on the component of the forecast along $r$.
\item If $a>0$, the expected SPO+ loss along $r$ is convex, and its minimizers are exactly the shifts $k$ that
satisfy the \emph{throughput-matching} condition
\begin{equation}\label{eq:matching-general}
\E\,T^-(2\hat u-p+2kv)\le\E\,T(z^\star)\le\E\,T^+(2\hat u-p+2kv),
\end{equation}
where
$T^-$ and $T^+$ are the smallest and largest values of $T$ among best decisions at the given forecast. If $\min_{\mathcal Z}T<\E\,T(z^\star)<\bar T$, the minimizers form a nonempty closed and bounded
interval. When the helper problem has larger expected $T$ than the perfect-foresight decision at $k=0$,
every minimizer is positive and equals $\delta^\star/(2a)$ for a matching cost increase $\delta^\star$ on
$T$. If $a=0$, the expected loss is constant along $r$.
\item With a squared-error term $\lambda\,\E[H^{-1}\|\hat u+kv-p\|^2]$, $\lambda>0$, on the unprocessed
forecast, the minimizer $k_\lambda$ is unique and lies within $a\bar TH/(\lambda\|r\|^2)$ of the
squared-error minimizer $-r^\top\E(\hat u-p)/\|r\|^2$.
\item Suppose $g$ saturates along $r$ at $p$: $g(p+kv)=s$ for all $k\ge k_0$. If $a=0$, then
$\ell(p+kv,p)=0$ for every $k$, while for $k\ge k_0$ the deployed regret is that of a best decision at the
fixed vector $s$, whatever $k$ is; for the battery with $s=\bar p\1$ this is the whole of $\Pi^\star(p)$.
\item Suppose $g$ acts coordinate by coordinate, $g(u)_t=g_t(u_t)$, and let $S$ be a set of
coordinates such that each $g_t$, $t\in S$, is constant on an open interval containing $u_t$. Then
$u\mapsto\ell(g(u),p)$ does not change when $(u_t)_{t\in S}$ move within these intervals, so
every chain-rule subgradient of it has zero entries in $S$. Clipping to a box satisfies this, with $S$ the
coordinates strictly outside the bounds. For maps that mix coordinates, such as projection onto a simplex,
saturation of one output coordinate does not imply this. A sufficient local condition for a general map $g$
that is differentiable at $u$ is that the whole $j$-th column of its Jacobian vanishes, $Dg(u)e_j=0$;
then every chain-rule subgradient of $\ell(g(\cdot),p)$ at $u$ has a zero $j$-th entry.
\item For any $w$ with $w^\top r\ne0$ and any level estimate $\hat L$, the anchored forecast
$u-(w^\top u/w^\top r)\,r+\hat L\,r$ does not change when $u$ moves along $r$, so training on it
cannot learn a shift along $r$.
\end{enumerate}
\end{theorem}

Theorem~\ref{thm:general} rests on three assumptions: profit affine in the forecast; a feasible set that does not depend on the forecast and whose convex hull has finitely many extreme points; and the near-conservation identity of Assumption~\ref{ass:cons} with bounded $T$. Parts (i) and (iii), the convexity and matching statements in (ii), and part (vi) follow from these by the same arguments as the battery results below, with $\1$ replaced by $r$ and $(1-\eta)T$ by $aT$; the constant $b$ cancels (Section~\ref{app:proofs}). The existence claim in (ii) needs, in addition, that $T$ attains its minimum and maximum on $\mathcal V_{\mathcal Z}$, which holds because $\mathcal V_{\mathcal Z}$ is finite, and a dominated-convergence step that uses the bound $0\le T\le\bar T$. Parts (iv) and (v) concern the post-processing map and need only the stated properties of $g$. The mechanism is therefore not specific to price caps or to batteries. The theorem says where the problem comes from: a direction in
which the decision loss is flat or nearly flat ($a$ zero or small) combined with a post-processing map that
is not neutral along that direction. Two examples besides batteries:
\begin{itemize}[leftmargin=1.2em]
\item \emph{Batteries.} By Lemma~\ref{lem:cons}, a battery that ends the day at its initial charge satisfies
Assumption~\ref{ass:cons} with $r=\1$, $a=1-\eta$, $T$ its charged energy and $b=0$. Sections~\ref{sec:noloss}
to~\ref{sec:training} develop this case.
\item \emph{Flexible load.} A process or a fleet of electric vehicles must buy a fixed amount of energy $E$ over the day, $\sum_t\Delta x_t=E$ with $0\le x_t\le\bar x_t$, at forecast prices. With $r=\1$ and $q(z)=-x$, Assumption~\ref{ass:cons} holds with $a=0$, so decision-focused losses carry no information at all about the price level. A two-period example: one unit must be bought in either period, actual prices are $(5,10)$ and the cap is 10. The forecast $(12,17)$, which is the actual price vector shifted up by 7, has zero SPO+ loss by part (iv) and leads to the correct purchase in period 1. After clipping it becomes $(10,10)$, at which both periods are optimal; a solver that picks period 2 pays 10 instead of 5. If every block is clipped, every feasible schedule is optimal at the clipped forecast, and the deployed schedule is decided by tie-breaking alone.
\item \emph{Storage with proportional losses.} Inventory or a water reservoir with a fixed terminal level and
proportional losses satisfies the assumption with $a$ equal to the loss rate, as the battery does with
$a=1-\eta$.
\end{itemize}
We test only the battery case.

\subsection{Batteries with no energy loss}\label{sec:noloss}
\begin{theorem}[Level invariance and failure after projection]\label{thm:eta1}
Let $\bar T=\max_{z\in\F}T(z)$.
\begin{enumerate}[label=(\roman*),leftmargin=2em]
\item For all $\hat p,p$ and $k$,
\begin{equation}\label{eq:levelbound}
\big|\ell(\hat p+k\1,p)-\ell(\hat p,p)\big|\le2|k|(1-\eta)\bar T.
\end{equation}
In particular, $\ell$ does not change with the level when $\eta=1$.
\item Suppose $\omega_b,\omega_s>0$, the cap $\bar p>0$ is the same in all blocks of the day,
$p\in[0,\bar p]^H$ and $\Pi^\star(p)>0$. For every $k\ge\bar p-\min_tp_t$, the projected forecast
$\Pi_{\mathcal B}(p+k\1)=\bar p\1$ has idling as its only best schedule. So
$\R^{\mathrm{clip}}(p,p+k\1)=\Pi^\star(p)$, while $\ell(p+k\1,p)\le2k(1-\eta)\bar T$, which is
zero when $\eta=1$.
\item Take a model class that allows level shifts (for example, any model with an intercept), and let
$\eta=1$. Then the training SPO+ loss is the same for all level shifts. So some of the models with the lowest training loss give projected forecasts that lose the whole perfect-foresight profit on every
training day.
\end{enumerate}
\end{theorem}

Part (ii) shows the worst case. Its regret statement holds for every $\eta$: if every block of the
shifted forecast is at or above the cap, the battery does not trade and loses its whole profit. Only the
statement that the unregularized SPO+ loss is zero needs $\eta=1$; it does not cover the regularized
objective used in training. As a result, no bound of the form $\R^{\mathrm{clip}}(p,\hat p)\le\phi(\ell(\hat p,p))$ with $\phi(0)=0$ can hold for the raw SPO+ loss once forecasts are projected: with $\eta=1$, part (ii) gives a forecast with zero raw loss whose regret after projection is the whole perfect-foresight profit. A surrogate guarantee for the raw decision problem does not transfer through the deployment map; only the loss evaluated at the deployed, projected forecast bounds the deployed regret. For real batteries ($\eta<1$), part (i) only bounds how much the loss can change, and does not by itself predict a failure.
The next result explains what happens when $\eta<1$.

\subsection{Batteries with energy loss: an upward level shift}
For a price vector $\pi$, let $T^-(\pi)$ and $T^+(\pi)$ be the smallest and largest charged energy among
the best schedules at $\pi$. They are equal when the best schedule is unique.

\begin{lemma}[Charging falls as prices rise]\label{lem:mono}
If $\eta<1$, $k_1<k_2$ and $z_i$ is a best schedule at $\pi+k_i\1$, then $T(z_1)\ge T(z_2)$.
\end{lemma}

In all statements about expected losses (Theorem~\ref{thm:bias},
Corollary~\ref{cor:reg}, Proposition~\ref{prop:decomp} and Theorem~\ref{thm:general}) we assume that prices and
forecasts have finite first moments, $\E\|p\|_1<\infty$ and $\E\|\hat p\|_1<\infty$, and, where a squared-error
term is used, finite second moments. Because the feasible set has finitely many extreme points,
$|\Pi^\star(\pi)|\le c_0+c_1\|\pi\|_1$ for constants $c_0,c_1$, so all losses are integrable, and charged energy is
bounded; these facts justify exchanging expectation with one-sided differentiation and limits in the proofs.

\begin{assumption}[The helper problem has excess throughput]\label{ass:over}
At the forecast's own level, the helper problem charges more energy on average than the perfect-foresight
schedule:
\begin{equation}\label{eq:over}
\E\,T^-(2\hat p-p)>\E\,T(z^\star).
\end{equation}

\end{assumption}

Assumption~\ref{ass:over} is a condition, not a general property of SPO+. It holds in both of our
markets (Section~\ref{sec:results} and Section~\ref{app:de}), but it can fail, for example for a forecast
that is too smooth. Without it, Theorem~\ref{thm:bias} still gives the matching condition, but the best
shift can be zero or negative.

\begin{theorem}[Loss-minimizing level shift of the SPO+ loss]\label{thm:bias}
Let $\eta<1$. Fix a forecast rule $\hat p$ that uses only the information available before bidding, and
let $\Phi(k)=\E\,\ell(\hat p+k\1,p)$ be the expected SPO+ loss when the forecast level is shifted by $k$.
Let $K^\star=\argmin_k\Phi(k)$.
\begin{enumerate}[label=(\roman*),leftmargin=2em]
\item $\Phi$ is convex, with one-sided derivatives
\begin{equation}\label{eq:phideriv}
\Phi'_\pm(k)=-2(1-\eta)\,\E\big[T^{\mp}(2\hat p-p+2k\1)-T(z^\star)\big].
\end{equation}

\item \emph{Throughput matching.} $k\in K^\star$ if and only if
\begin{equation}\label{eq:matching}
\E\,T^-(2\hat p-p+2k\1)\le\E\,T(z^\star)\le\E\,T^+(2\hat p-p+2k\1).
\end{equation}
The expected charged energy of the
helper problem does not increase as $k$ increases. If $0<\E\,T(z^\star)<\bar T$, then $K^\star$ is a
nonempty, closed and bounded interval.
\item \emph{Direction.} Under Assumption~\ref{ass:over}, every $k\in K^\star$ is positive.
\item \emph{Size.} For $\delta\in\mathbb R$, let $T^\pm_\delta$ be the charged energy (largest and
smallest over the best schedules) when the buying cost is raised by $\delta$. Let
\[
D^\star=\big\{\delta:\ \E\,T^-_\delta(2\hat p-p)\le\E\,T(z^\star)\le\E\,T^+_\delta(2\hat p-p)\big\}
\]
be the set of matching cost increases. Then $K^\star=\{\delta/(2(1-\eta)):\delta\in D^\star\}$. If
$\min D^\star\ge\delta_0>0$ for all $\eta$ near one, then $\min K^\star\to\infty$ as $\eta\to1$.
\item \emph{Fenchel--Young.} For $L_\sigma$, (i)--(iv) hold with comparison price $\hat p+\sigma Z$ (and
the expectation also taken over $Z$), and with $2(1-\eta)$ replaced by $1-\eta$.
\end{enumerate}
\end{theorem}

\begin{remark}[Why the helper problem has excess throughput]\label{rem:guard}
The comparison price $2\hat p-p=\hat p+(\hat p-p)$ contains the forecast error twice. Errors can make the
price differences within the day look larger, so the helper problem may find more profitable cycles than actually exist. The perturbed loss adds noise to the forecast and can have the same effect. Whether this
happens depends on the forecast; Assumption~\ref{ass:over} states it as a condition. By
Lemma~\ref{lem:cons}, raising the level is the same as raising the cost of charging. So DFL uses the
forecast level to reduce trading, which works like a safety margin \citep{matsumoto2026guard}. Without a
binding cap this acts as a small trading cost. With a binding cap it pushes the forecast peaks above the
cap.
\end{remark}

Practical DFL methods, including all five we test, add a
squared-error term to the decision loss. The next result covers the case where both terms use the
forecast \emph{before} projection, as in V1. The case where both terms use the forecast after projection,
as in V2--V5, is different and is discussed in Remark~\ref{rem:projpen}.

\begin{corollary}[SPO+ with a squared-error term, before projection]\label{cor:reg}
Let $\lambda>0$ and let $\bar e=\E\big[H^{-1}\1^\top(\hat p-p)\big]$ be the mean forecast error. Let
\[
\mathcal L_\lambda(k)=\Phi(k)+\lambda\,\E\Big[\tfrac1H\big\|\hat p+k\1-p\big\|^2\Big],
\]
with both terms in the same money units. Then:
\begin{enumerate}[label=(\roman*),leftmargin=2em]
\item $\mathcal L_\lambda$ is strictly convex and has a unique minimizer $k_\lambda$, given by
\[
(1-\eta)\big(\E\,T^-(2\hat p-p+2k_\lambda\1)-\E\,T(z^\star)\big)\ \le\ \lambda\,(k_\lambda+\bar e)\ \le\
(1-\eta)\big(\E\,T^+(2\hat p-p+2k_\lambda\1)-\E\,T(z^\star)\big).
\]
\item $|k_\lambda+\bar e|\le(1-\eta)\bar T/\lambda$. So the shift away from the squared-error level
$-\bar e$ stays bounded, and it goes to zero as $\eta\to1$ or $\lambda\to\infty$. Unlike the minimizer of
the pure loss, which can diverge under the condition of Theorem~\ref{thm:bias}(iv), it stays bounded as
$\eta\to1$.
\item If the helper problem has excess throughput at the squared-error level, that is,
$\E\,T^-(2\hat p-p-2\bar e\1)>\E\,T(z^\star)$, then $k_\lambda>-\bar e$: the shift is upward. If, in
addition, $-\bar e<\min K^\star$, then $-\bar e<k_\lambda\le\min K^\star$.
\end{enumerate}
\end{corollary}

\begin{remark}[Squared error after projection: no guarantee]\label{rem:projpen}
In V2--V5, the decision loss and the squared-error term both use the projected forecast:
\[
\mathcal L^{\mathrm{proj}}_\lambda(k)=\E\Big[\ell\big(\Pi_{\mathcal B}(\hat p+k\1),p\big)
+\tfrac{\lambda}{H}\big\|\Pi_{\mathcal B}(\hat p+k\1)-p\big\|^2\Big].
\]
This function is generally not convex in $k$, and it is completely flat in $k$ once every block is clipped
(Proposition~\ref{prop:plateau}). So Corollary~\ref{cor:reg} does not apply: there need not be a unique or a
bounded loss-minimizing level. The derivative of the penalty is

\begin{equation}\label{eq:projpen}
\frac{2\lambda}{H}\,\E\sum_t\chi_t(k)\big(\Pi_{\mathcal B}(\hat p+k\1)_t-p_t\big),\qquad \chi_t(k)=\mathbf 1\{\underline p<\hat p_t+k<\bar p\},
\end{equation}
so the penalty gives no pull back in clipped blocks.
\end{remark}

Theorem~\ref{thm:bias} and Corollary~\ref{cor:reg} describe the \emph{optimum} of the loss. They do not
say that a training algorithm reaches it. In a weakly curved direction, a stochastic optimizer with a
limited number of steps can stay far from the optimum. Section~\ref{sec:results} shows that this is what
happens in our experiments: with more training steps, the level error falls sharply.

\subsection{When projection preserves the decision}\label{sec:safe}
Let $u$ be a raw forecast, $v=\Pi_{\mathcal B}(u)$ its projection, and $\mathcal C=\{t:u_t\ne v_t\}$ the clipped
blocks. Let $z_u$ be a best schedule at $u$ and $z_v$ a best schedule at $v$. Projection changes prices
only in $\mathcal C$, and it lowers them in blocks cut at the cap and raises them in blocks raised to the floor.

\begin{proposition}[Projection and the decision]\label{prop:margin}
\begin{enumerate}[label=(\roman*),leftmargin=2em]
\item \emph{Exact condition.} $z_u$ is a best schedule at $v$ if and only if, for every $z\in\mathcal V$,
\[
J(u,z_u)-J(u,z)\ \ge\ \Delta\sum_{t\in\mathcal C}(u_t-v_t)\big(q_t(z_u)-q_t(z)\big);
\]
equivalently, $v$ lies in the normal cone
\begin{equation}\label{eq:normalcone}
N(z_u)=\big\{\pi:z_u\in\argmax_{z\in\F}J(\pi,z)\big\},
\end{equation}
which is a polyhedron.
\item \emph{Margin condition.} Suppose $z_u$ is the unique best point of $\mathcal V$ at $u$, and for every other $z\in\mathcal V$,
\[
J(u,z_u)-J(u,z)>\Delta\sum_{t\in\mathcal C}|u_t-v_t|\,|q_t(z_u)-q_t(z)|.
\]
Then $z_u$ is also the unique best schedule at $v$, so $\R^{\mathrm{clip}}(p,u)=\R^{\mathrm{raw}}(p,u)$.
\item \emph{Change in regret from projection.} For any choice of $z_u$ and $z_v$,
\[
\Delta(p-u)^\top\big(q(z_u)-q(z_v)\big)\ \le\ \R^{\mathrm{clip}}(p,u)-\R^{\mathrm{raw}}(p,u)\ \le\
\Delta(p-v)^\top\big(q(z_u)-q(z_v)\big).
\]
The two bounds differ by $\Delta\sum_{t\in\mathcal C}(u_t-v_t)\big(q_t(z_u)-q_t(z_v)\big)\ge0$, which involves only the
clipped blocks.
\end{enumerate}
\end{proposition}

Part (i) is exact: projection preserves the optimality of $z_u$ exactly when the change it makes to the
forecast keeps the forecast inside the normal cone of $z_u$. Parts (i) and (ii) are standard sensitivity
arguments. None of the three parts is specific to batteries. They hold for any decision problem whose profit is affine in the forecast and whose feasible set has a convex hull with finitely many extreme points, including mixed-integer sets such as our physical $\F$. They also hold for any post-processing map, with $\mathcal C$ the set of coordinates it changes. Checking (i) or (ii) directly, however, requires a comparison with every
point of $\mathcal V$, or a new solve at $v$. Part (iii) bounds the change in regret caused by projection, which can be positive, negative or zero,
because projection can improve as well as worsen a deployed decision. Both bounds depend on forecast errors
in all blocks; only the gap between them is confined to the clipped blocks. The next result uses the structure of the battery
problem to give a sufficient condition that can be checked block by block from the raw solution alone.

Suppose every price in $u$ and $v$ exceeds $\theta$, so the LP is
exact (Lemma~\ref{lem:lp}). Let $\mu$ be an optimal dual solution of the LP at $u$, with $\mu_t$, the value of stored energy in block $t$, the dual variable of the storage equation of block $t$. When the dual solution is not unique, the thresholds below depend on this choice. Define the \emph{charge threshold}
$\theta^c_t=\eta_c\mu_t-\omega_b$ and the \emph{discharge threshold} $\theta^d_t=\mu_t/\eta_d+\omega_s$. By
the optimality conditions, at $u$ the battery charges at full power in blocks where $u_t<\theta^c_t$,
discharges at full power where $u_t>\theta^d_t$, is idle where $\theta^c_t<u_t<\theta^d_t$, and runs at partial
power only where $u_t$ equals a threshold. The value $\mu_t$ changes only when the storage is full or empty.

\begin{theorem}[Safe clipping for batteries]\label{thm:safe}
Suppose every price in $u$ and $v$ exceeds $\theta$, and that, for some optimal dual solution $\mu$ of the LP at $u$, every clipped block $t\in\mathcal C$ satisfies:
\begin{enumerate}[label=(\alph*),leftmargin=2em]
\item $z_u$ does not charge or discharge at partial power in block $t$;
\item if block $t$ is cut to the cap: when $z_u$ discharges at full power there, $\bar p\ge\theta^d_t$; when
$z_u$ is idle there, $\bar p\ge\theta^c_t$;
\item if block $t$ is raised to the floor: when $z_u$ charges at full power there, $\underline p\le
\theta^c_t$; when $z_u$ is idle there, $\underline p\le\theta^d_t$.
\end{enumerate}
Then $z_u$ is a best schedule at $v$. If it is the only best schedule at $v$ (or the selection rule picks
it), then $\R^{\mathrm{clip}}(p,u)=\R^{\mathrm{raw}}(p,u)$. Otherwise the solver may return a different schedule
that is equally good at $v$ but not at the actual prices.
\end{theorem}

Because any optimal primal solution and any optimal dual solution of an LP satisfy complementary
slackness together, the certificate may use any optimal dual solution, and it holds if it holds for at least
one. Our computations use the dual solution returned by the solver, which makes the computed certificate
conservative (Section~\ref{sec:results}). In words: after a forecast peak is cut to the cap, the raw schedule remains optimal as long as the cap is still above
the price at which discharging in that block is worthwhile.

\begin{corollary}[An unclipped partial-power block protects the peak]\label{cor:marginal}
Suppose $z_u$ discharges at partial power in an unclipped block $m$, and the storage stays strictly between
empty and full from block $m$ to a clipped block $t$ in which $z_u$ discharges at full power. Then condition
(b) holds in block $t$: $\theta^d_t=\theta^d_m=u_m\le\bar p$.
\end{corollary}

So, after the full-power discharge blocks are clipped, the raw schedule remains optimal as long as the block where the battery runs at partial
power, which sets the discharge threshold, is not clipped. A 1\,MW/4\,MWh battery discharges in about 16
to 23 of the 96 blocks on a typical day in our data. Once more blocks than this are clipped, the clipped
set must include blocks that matter at the margin. Section~\ref{sec:results} checks Theorem~\ref{thm:safe}
on the trained forecasts and links it to how often clipping is costly.

Theorem~\ref{thm:safe} guarantees that the raw schedule stays
optimal after projection, not that a solver returns it. A deterministic selection rule closes this gap.

\begin{corollary}[Schedule-preserving selection]\label{cor:select}
Consider the rule that, after projection, deploys $z_u$ whenever $z_u$ is a best schedule at $v$, and
otherwise deploys a best schedule at $v$. Under this rule, the deployed regret equals $\R^{\mathrm{raw}}(p,u)$
whenever $z_u$ is a best schedule at $v$, in particular whenever the conditions of Theorem~\ref{thm:safe} hold.
Otherwise, its difference from $\R^{\mathrm{raw}}(p,u)$ satisfies the bounds of Proposition~\ref{prop:margin}(iii).
\end{corollary}

The rule needs only the raw schedule and one check of its value at the projected forecast; the certificate
of Theorem~\ref{thm:safe} gives a sufficient condition that can be read off the raw solution alone. Under this
rule, projection can change the deployed regret only on days when the raw schedule is no longer optimal
after projection. On those days, a natural refinement chooses, among the schedules that are optimal at
$v$, one that maximizes profit at the raw forecast $u$ (a lexicographic rule). It uses only information
available before dispatch. It gives the same schedule as the rule of Corollary~\ref{cor:select} whenever $z_u$ is optimal at $v$, provided that ties among schedules that are optimal at both $v$ and $u$ are broken in favor of $z_u$. Otherwise the two rules can pick different schedules that are equally good at both forecasts. Our implementation breaks such ties in favor of $z_u$: when $z_u$ is optimal at $v$, it is
deployed directly. Section~\ref{sec:rule} evaluates both.

\subsection{Training on projected forecasts and anchoring}\label{sec:training}
An obvious idea is to train on the projected forecast, so that training matches how the forecast is
used. Then the training loss is still an upper bound on regret, because $\ell(\Pi_{\mathcal B}(u),p)\ge
\R^{\mathrm{clip}}(p,u)$ for every model output $u$. But this does not by itself stop the drift.

\begin{proposition}[Zero gradients after projection]\label{prop:plateau}
Let $\hat p=\Pi_{\mathcal B}(u)$ and let $\mathcal C(u)=\{t:u_t>\bar p\text{ or }u_t<\underline p\}$ be the
blocks that are cut. Then $\ell(\Pi_{\mathcal B}(u),p)$ does not change with small changes in
$(u_t)_{t\in\mathcal C(u)}$, so the gradient is zero in those blocks. If $u_t\ge\bar p$ for all $t$, the
loss is the same for all $k\ge0$ along $u+k\1$, and equals the loss of the flat forecast $\bar p\1$.
\end{proposition}

Proposition~\ref{prop:plateau} is about the \emph{direct} gradient of the clipped blocks. A regression or a
neural network shares parameters across blocks. Gradients from unclipped blocks can therefore still change the
parameters and move clipped forecasts back inside the limits. But the loss gives no local signal from the
clipped blocks themselves. In our experiments, this indirect effect was not enough to prevent the drift
(Section~\ref{sec:results}).

\begin{corollary}[Anchored training on projected forecasts]\label{cor:anchor}
Let a model give raw outputs $u$. Form the anchored forecast
\begin{equation}\label{eq:anchor}
\tilde u=u-H^{-1}(\1^\top u)\,\1+\hat L\1,
\end{equation}
whose daily level $\hat L$ comes from a separate model, and deploy its projection $\hat p=\Pi_{\mathcal B}(\tilde u)$. Train on $\ell(\hat p,p)$ or $L_\sigma(\hat p,p)$, evaluated at the deployed forecast $\hat p$. Then:
\begin{enumerate}[label=(\roman*),leftmargin=2em]
\item the training loss does not change when $u$ is shifted up or down, so training cannot learn a common
level shift; this removes the freely learned level shift, but not every possible clipping error;
\item for SPO+, the training loss $\ell(\hat p,p)$, evaluated at the projected forecast $\hat p$ that is used for the decision, is an upper bound on the regret after projection on every training day;
\item the anchored forecast $\tilde u$ has daily mean exactly $\hat L$, so its level error comes only from the separate level model, not from DFL training; by Proposition~\ref{prop:bound}, it matters only through the factor $1-\eta$ unless projection cuts it. After projection, clipping can change the daily mean of $\hat p$.
\end{enumerate}
\end{corollary}

\begin{algorithm}[t]
\caption{Anchored decision-focused learning on projected forecasts (one refit)}\label{alg:anchor}
\begin{algorithmic}[1]
\Require training days $j$ with features, actual prices $p_j$, fixed perfect-foresight schedules
$z^\star_j$; a level model (e.g.\ OLS) giving $\hat L_j$; a forecasting model $u_\kappa$ with parameters $\kappa$; loss weight
$\lambda$
\For{each training step}
\State Sample a batch of days. For each day $j$:
\State \quad $u\gets u_\kappa(j)$;\quad $\tilde u\gets u-H^{-1}(\1^\top u)\1+\hat L_j\1$;\quad $\hat p\gets\Pi_{\mathcal B}(\tilde u)$
\State \quad $\tilde z\gets$ a best schedule at the comparison price $2\hat p-p_j$
\State \quad $\xi\gets2\Delta\big(q(\tilde z)-q(z^\star_j)\big)+\tfrac{2\lambda}{H}(\hat p-p_j)$ \Comment{gradient of $\ell+\tfrac{\lambda}{H}\|\hat p-p_j\|^2$}
\State \quad $\xi\gets \xi\odot\mathbf 1\{\underline p<\tilde u<\bar p\}$;\quad $\xi\gets \xi-H^{-1}(\1^\top \xi)\1$ \Comment{projection, then centering}
\State \quad add $\nabla_\kappa u_\kappa(j)^\top \xi$ to the total gradient
\State Update $\kappa$ (Adam).
\EndFor
\State Use the forecast $\hat p_D=\Pi_{\mathcal B}\big(u_\kappa(D)-H^{-1}(\1^\top u_\kappa(D))\1+\hat L_D\1\big)$
\end{algorithmic}
\end{algorithm}

\section{Data and test design}\label{sec:design}
We now turn from theory to data. This section describes the market data and the battery, the forecasting
methods we compare, and the test design and statistics used in Section~\ref{sec:results}.

\subsection{Market, data and battery}
We use 15-minute prices from the IEX day-ahead market \citep{iexDAMMarketData}, from 1~April~2022 to
31~August~2026 for the main study (1{,}614 calendar days); Section~\ref{sec:fresh} uses later data. The
regulator's order of 1~April~2022 set a ceiling of Rs\,12/kWh (Rs\,12{,}000/MWh) for the day-ahead and
real-time markets \citep{cerc2022sm4}, extended to further market segments on 6~May~2022 \citep{cerc2022sm5}.
The April order does not state the first delivery day it affects. In our data, prices on 1 and 2~April~2022
still exceed Rs\,12{,}000/MWh (up to Rs\,20{,}000/MWh), but no later price does, so our panel starts on
3~April~2022. The ceiling was lowered to Rs\,10/kWh for delivery from 4~April~2023 \citep{cerc2023sm4}, and
the floor is zero. We use the ordinary day-ahead market, not its high-price segment, which has a separate
ceiling of Rs\,20/kWh. After a 90-day warm-up, the remaining 1{,}612 days give \IEXnDays{} delivery days, from
2~July~2022 to 31~August~2026. A day $D$ is evaluated only if days $D-7$, $D-6$, $D-1$, $D$ and $D+1$ all have
complete data. This rule comes from the original design, in which every day also had to support a 48-hour
schedule that uses day $D+1$, and we keep it unchanged; all 1{,}522 days satisfy it.

The battery is 1\,MW/4\,MWh with $\eta=0.85$ ($\eta_c=\eta_d=\sqrt{0.85}$), a wear cost of Rs\,1{,}062.5 per
MWh cycled, split between buying and selling, an exchange fee of Rs\,20/MWh for each
\citep{cerc2023transactionfee}, and $\bar x=0$. With perfect foresight (the ``oracle''), it would earn
Rs\,\IEXoracle\ lakh per MW-year. Section~\ref{app:de} reports a second test on the German day-ahead market.

\subsection{Forecasting methods}
All methods use only prices up to day $D-1$. Every forecast is projected to the price limits and then
used in the same 24-hour dispatch optimization.
\begin{itemize}[leftmargin=1.2em]
\item \emph{Seasonal naive}: the prices of day $D-1$.
\item \emph{Regression}: $p_{D,t}=\beta_0+\beta_1p_{D-1,t}+\beta_2p_{D-7,t}+\text{block}_t
+\alpha_1\sin(2\pi m/12)+\alpha_2\cos(2\pi m/12)+\varepsilon$, estimated by least squares on all past data.
\item \emph{Gradient boosting}: histogram gradient boosting with settings fixed in advance. Inputs are
lagged prices, block, day of week and month. Refit every week.
\item \emph{LEAR-type LASSO} \citep{lago2021forecasting,tibshirani1996lasso}: one LASSO model for each
block, on asinh-transformed prices. Inputs are all prices of day $D-1$; the same block on days $D-2$,
$D-3$ and $D-7$; daily summary values; and day-of-week dummies. It uses the last 728 days, chooses
variables by AIC, and is refit every week.
\item \emph{EW3}: the simple average of the projected regression, gradient boosting and LASSO forecasts
\citep{bates1969combination}.
\item \emph{DFL methods}. All are trained on the decision loss plus $\lambda$ times the squared error of
the forecast, using Adam \citep{kingma2015adam}. They start from the regression (the neural network is
first trained on squared error) and are refit every 28 days on all past data, using batches of 16 days.
\begin{itemize}[leftmargin=1.2em]
\item V1: SPO+ with the regression's inputs, trained on the forecast before projection.
\item V2: as V1, but trained on the forecast after projection.
\item V3: as V2, plus anchoring to the regression's daily level (Algorithm~\ref{alg:anchor}).
\item V4: a neural network (MLP, one hidden layer of 32 units) trained with SPO+, with anchoring and
projection.
\item V5: a perturbed Fenchel--Young method with anchoring, using 4 random draws.
\end{itemize}
\end{itemize}

\subsection{Test design and statistics}\label{sec:testdesign}
\textbf{Test design.} We fixed the following before the tuning and holdout runs of the DFL models (see the last paragraph of this section).
\begin{itemize}[leftmargin=1.2em]
\item \emph{Periods.} The \emph{validation year} is 2~July~2022 to 30~June~2023 (364 days). It is used
only for tuning. The \emph{holdout} is 1~July~2023 to 31~August~2026 (1{,}158 days). It is used only for
testing.
\item \emph{Tuning.} Each DFL method is tuned on its own grid, using only validation-year regret:
$\lambda\in\{0.1,1\}$ and learning rate in $\{3\times10^{-4},10^{-3},3\times10^{-3}\}$. For V5, we use
$\lambda=0.1$, $\sigma\in\{300,1000\}$ Rs/MWh and the same learning rates.
\item \emph{Seeds.} The chosen setting is run with three random seeds. The statistical methods have no
settings to tune.
\item \emph{Hindsight check.} As an extra check, we also report the holdout regret of every setting for the two best DFL methods, as if the settings had been chosen with hindsight. This favors DFL, so it is a strict test of our conclusion.
\end{itemize}

On 39\% of days there is more than one best schedule, mostly because many
blocks clear at exactly the cap. The SPO+ loss depends on which one is used, so we fix the choice for
each day.

\textbf{Statistics.} We compare methods by their daily regret differences. Confidence intervals come from
a circular block bootstrap with 7-day blocks and 4{,}000 repetitions
\citep{kunsch1989bootstrap,politis1992circular}. All holdout days are under the Rs\,10 cap. Regret is
given in lakh rupees per MW-year.

The LASSO and EW3 benchmarks were first found
in earlier exploratory work on the full Indian sample. Earlier versions of this study, and our own
exploratory runs of a linear SPO+ variant on 23~October~2026, also used the full evaluation range, including
the later holdout. The holdout is therefore not a fully independent test of these choices. The DFL benchmark
protocol was fixed in code later on 23~October~2026: the variants, tuning grids and selection rule before the
tuning runs, and the seeds of the final runs before any holdout run. The German test was specified on
24~October~2026, before any German result; it checks the ranking in another market, but it cannot test the
cost of projection, because German prices never reach their limits. The fresh-test protocol
(Section~\ref{sec:fresh}) was written on 25~October~2026. Its days had already occurred, and we had checked
the new data for completeness and computed the share of blocks at the cap, but no forecast or regret for
those days had been computed; the protocol was fixed before the analysis, not registered in advance. The
operating rules of Section~\ref{sec:rule} were fixed on 2~November~2026, each before its own comparison, but
after the forecasting results were known, so their results are exploratory. For the fresh test and the
operating rules, the original time-stamped protocol files are in the replication material. For the DFL
benchmark and the German test, the design is recorded in the scripts that implement it, whose modification
times precede the results. These times document the order of the work but are not independently certified.

\section{Empirical results}\label{sec:results}
This section first examines the predictions of Section~\ref{sec:theory} in the Indian data
(Section~\ref{sec:theorydata}). It then compares all forecasting methods on the holdout
(Sections~\ref{sec:main} and~\ref{sec:designs}), reports a fresh test on later data (Section~\ref{sec:fresh}),
evaluates operating rules based on the certificate of Theorem~\ref{thm:safe} (Section~\ref{sec:rule}), and
ends with the size of the gains (Section~\ref{sec:gains}).

\subsection{Level shift and projection in the data}\label{sec:theorydata}
At the regression forecast, the SPO+ helper problem
charges 5.25~MWh per day, while the oracle charges 5.04. So Assumption~\ref{ass:over} holds, and
Theorem~\ref{thm:bias} predicts an upward loss-minimizing shift.

\begin{figure}[!htbp]\centering
\includegraphics[width=0.9\linewidth]{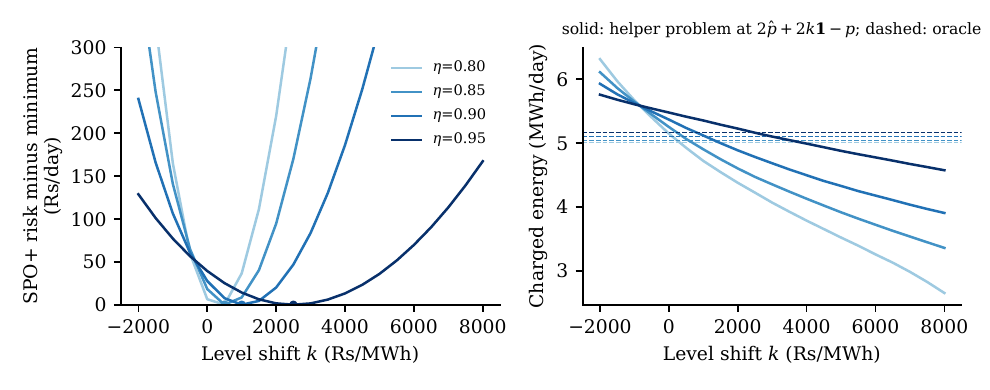}
\caption{Left: expected pure SPO+ loss for level shifts of the regression forecast, on a 500 Rs/MWh grid
(lowest grid point marked; exact minimizers in Table~\ref{tab:eta}). Right: charged energy of the SPO+
helper problem (solid) and of the oracle (dashed). The lines cross at the loss-minimizing level shift of
Theorem~\ref{thm:bias}.}\label{fig:eta}
\end{figure}

For the regression forecast, Table~\ref{tab:eta} gives the minimizer of the
expected SPO+ loss $\Phi(k)$ for four efficiencies, each with its own oracle and wear-cost split, and
Figure~\ref{fig:eta} plots $\Phi$ and the helper problem's charged energy against the shift. We compute the
minimizers by bisection on the condition in Theorem~\ref{thm:bias}(ii), with the one-sided charged energies
$T^\mp$ obtained by raising or lowering the buying cost by 0.01 Rs/MWh, to a tolerance of 1 Rs/MWh. The lower
and upper ends of $K^\star$ agree to within this tolerance, so the minimizer is effectively a single point. It
is positive and rises from \etaKa\ to \etaKd\ Rs/MWh as $\eta$ goes from 0.80 to 0.95, while the matching cost
increase $\delta^\star=2(1-\eta)k^\star$ moves only from 130 to 250. The pure loss is almost flat in this
direction: moving from $k=0$ to $k^\star$ lowers it by only 9 to 39 Rs per day, out of about 5{,}500. With
the squared-error term of Corollary~\ref{cor:reg}, as in V1 (in our units the tuned methods use $\lambda=0.1$
or $\lambda=1$), the loss-minimizing shift is much smaller: 152 to 235 Rs/MWh ($\lambda=0.1$) or 16 to 34
Rs/MWh ($\lambda=1$) above the squared-error level, and it does not grow as $\eta\to1$. This result does not
apply to V2--V5 (Remark~\ref{rem:projpen}). Because the mean error of $\hat p+k\1$ is $\bar e+k$, these shifts
are on the same scale as the level errors in Tables~\ref{tab:drift} and~\ref{tab:conv}, although they refer
to a different forecast and sample.

\begin{itemize}[leftmargin=1.2em]
\item At every efficiency, the lower and upper ends of $K^\star$ agree to within this accuracy, so the
minimizer is a single point for practical purposes. It is positive and rises from \etaKa\ to \etaKd\
Rs/MWh as $\eta$ goes from 0.80 to 0.95.
\item The matching cost increase $\delta^\star=2(1-\eta)k^\star$ moves only from 130 to 250, while
$1/(1-\eta)$ grows four times.
\item Moving from $k=0$ to $k^\star$ lowers the loss by only 9 to 39 Rs per day, out of about 5{,}500.
The pure loss is almost flat in the level direction.
\item For the objective of Corollary~\ref{cor:reg} (SPO+ plus squared error, both before projection, as in
V1), the loss-minimizing shift is much smaller. In our units (thousand rupees, and thousand Rs/MWh for prices), the
tuned methods use $\lambda=0.1$ or $\lambda=1$ (V1 uses $\lambda=1$). Then the loss-minimizing shift is 152 to 235
Rs/MWh ($\lambda=0.1$) or 16 to 34 Rs/MWh ($\lambda=1$) above the squared-error level. It does not grow as
$\eta\to1$. This result does not apply to V2--V5 (Remark~\ref{rem:projpen}).
\item \emph{A common scale.} The mean error of the shifted forecast $\hat p+k\1$ relative to actual prices
is $\bar e+k$. So the shift of $k_\lambda$ above $-\bar e$ is exactly the mean forecast error relative to
actual prices at the optimum. This is the same scale as the level errors in Tables~\ref{tab:drift}
and~\ref{tab:conv}. The numbers still refer to different forecasts (the regression here, V1 there) and
different samples.
\end{itemize}

\begin{table}[!htbp]\centering\small
\begin{threeparttable}
\caption{Loss-minimizing level shift for the SPO+ loss, with and without the squared-error term. Regression
forecast, India (\IEXnDays{} days).}\label{tab:eta}
\footnotesize\setlength{\tabcolsep}{4pt}
\begin{tabular}{cccccccc}\toprule
& \multicolumn{2}{c}{Charged energy at $k=0$} & \multicolumn{3}{c}{Pure SPO+ loss} & \multicolumn{2}{c}{Shift of $k_\lambda$ above $-\bar e$}\\
& \multicolumn{2}{c}{(MWh/day)} & \multicolumn{3}{c}{(Theorem~\ref{thm:bias})} & \multicolumn{2}{c}{(Rs/MWh, Corollary~\ref{cor:reg})}\\
\cmidrule(lr){2-3}\cmidrule(lr){4-6}\cmidrule(lr){7-8}
$\eta$ & helper problem & oracle & $k^\star$ (Rs/MWh) & $\delta^\star$ (Rs/MWh) & \makecell{$\Phi(0)-\Phi(k^\star)$\\(Rs/day)} & $\lambda=0.1$ & $\lambda=1$\\\midrule
0.80 & 5.15 & 5.00 & 324 & 130 & 9.1 & 189 & 32\\
0.85 & 5.25 & 5.04 & 573 & 172 & 18.5 & 235 & 34\\
0.90 & 5.37 & 5.10 & 1{,}064 & 213 & 27.7 & 227 & 28\\
0.95 & 5.48 & 5.16 & 2{,}496 & 250 & 39.3 & 152 & 16\\
\bottomrule\end{tabular}
\begin{tablenotes}\footnotesize
\item Minimizers found by bisection on the one-sided matching conditions, with one-sided charged energies
computed at buying-cost changes of $\pm0.01$ Rs/MWh. Reported values are accurate to within
$\pm1$ Rs/MWh. At each $\eta$, the left and right ends of $K^\star$ agree within this accuracy. The mean
error of the regression forecast is $\bar e=79$ Rs/MWh. At $k=0$, $\E\,T^-$ and $\E\,T^+$ agree to two
decimals. Stability check at $\eta=0.85$: with perturbations of 0.001, 0.01 and 0.1 Rs/MWh, $k^\star=572.8$
in all three cases; with 1 Rs/MWh it is 569.8. LPs are solved with HiGHS at its default primal and dual
feasibility tolerances ($10^{-7}$).
\end{tablenotes}
\end{threeparttable}
\end{table}

Table~\ref{tab:drift} retrains the tuned V1 ($\lambda=1$, 150 steps
per refit) for three efficiencies. Its learned level moves up by about 2{,}400 to 2{,}500 Rs/MWh relative to
actual prices in the validation year, far more than the loss-minimizing 16 to 34 Rs/MWh (Table~\ref{tab:eta}).

Table~\ref{tab:conv} retrains V1 (seed 0) over the
validation year with 150 to 5{,}000 steps per refit. At each refit it also applies an \emph{exact level
correction}: using only that refit's training days, the level shift that minimizes the model's own training
objective (by bisection on the condition in Corollary~\ref{cor:reg}), applied until the next refit. From 150 to
5{,}000 steps, the out-of-sample level error falls from \ConvLvlA\ to \ConvLvlD\ Rs/MWh and the regret from
\ConvRegA\ to \ConvRegD\ lakh per MW-year, at a cost of about 35 thousand to one million training solves (half a
minute to under seven minutes on one CPU core). Most of the drift in our main runs therefore comes from
incomplete training. The level correction helps much less, because most of the level error appears only out
of sample. In-sample level errors are at most about 520 Rs/MWh with 150 steps, but with few steps the fitted
intercept and slopes give biased forecasts on new days: in one run (seed 1), the least-squares intercept fell
from 5.9 to 1.2 thousand Rs/MWh between July and October 2022, while V1's moved only from 6.0 to 4.8. Even with
5{,}000 steps and the correction, V1's first-year regret (\ConvRegDcor) exceeds that of the regression (4.45)
and of EW3 (4.08), although it is below that of anchored SPO+ with 150 steps (4.86).

To check that the ranking is not an artifact of
under-training, we reran all five DFL variants (seed 0) with 2{,}000 steps per refit, 13 times the main
budget, under a protocol fixed before the runs (Table~\ref{tab:budget}). More training removes most of the
first-year drift of the unanchored variants: the validation level error falls from +2{,}449 to +606 Rs/MWh for
V1, reproducing Table~\ref{tab:conv}, and from +2{,}637 to +821 for V2, and first-year regret falls from 6.60
to 4.62 and from 6.95 to 4.93. Holdout regret, however, barely moves for V1 to V3 (by at most 0.01) and rises
for V4 (2.89 to 3.10) and V5 (3.49 to 3.66), apparently because longer training overfits the decision loss.
For V4, the most seed-sensitive variant, the rise also holds for seeds 1 and 2 (2.88 to 3.11 and 2.91 to
3.01), and the three-seed mean, 3.07, remains worse than EW3 by 0.65 [0.50, 0.80]. Every variant remains worse
than EW3 on the holdout, by at least \BudgetMinGap\ lakh per MW-year with every interval excluding zero, so
under-training explains the early drift but not the ranking.

\begin{table}[!htbp]\centering\small
\begin{threeparttable}
\caption{Level shift of unconstrained SPO+ (V1) in actual training, India.}\label{tab:drift}
\begin{tabular}{cccccccc}\toprule
& \multicolumn{2}{c}{V1 level error (Rs/MWh)} & Raw blocks & \multicolumn{2}{c}{Validation regret} & \multicolumn{2}{c}{Holdout regret}\\
\cmidrule(lr){2-3}\cmidrule(lr){5-6}\cmidrule(lr){7-8}
$\eta$ & validation & holdout & above cap (val.) & V1 & Regression & V1 & Regression\\\midrule
0.80 & 2{,}433 & 88 & 8.7\% & 7.63 & 4.64 & 3.20 & 3.27\\
0.85 & 2{,}444 & 83 & 8.4\% & 6.52 & 4.45 & 3.09 & 3.17\\
0.95 & 2{,}466 & 45 & 8.3\% & 5.17 & 4.25 & 2.98 & 2.92\\
\bottomrule\end{tabular}

\begin{tablenotes}\footnotesize
\item V1 is retrained for each efficiency with the setting chosen at $\eta=0.85$. The $\eta=0.85$ row is
the average of two seeds. ``Level error'' is the average over days of (daily mean of the forecast before
projection) minus (daily mean of the actual price). Regret is in Rs lakh per MW-year.
\end{tablenotes}
\end{threeparttable}
\end{table}

\begin{table}[!htbp]\centering\small
\begin{threeparttable}
\caption{V1 with more training steps, validation year (seed 0).}\label{tab:conv}
\footnotesize\setlength{\tabcolsep}{2.6pt}
\begin{tabular}{lccccccccc}\toprule
& & & & & \multicolumn{3}{c}{Forecast as trained} & \multicolumn{2}{c}{With exact level correction}\\
\cmidrule(lr){6-8}\cmidrule(lr){9-10}
Steps per refit & \makecell{Training solves\\(thousands)} & \makecell{Time\\(s)} & \makecell{Training\\objective} & \makecell{In-sample\\level error} & \makecell{Out-of-sample\\level error} & \makecell{Blocks\\above cap} & Regret & \makecell{Out-of-sample\\level error} & Regret\\\midrule
150 & 35 & 32 & 10.23 & +311 & +2{,}449 & 8.5\% & 6.60 & +2{,}183 & 6.07\\
500 & 104 & 54 & 9.20 & +125 & +1{,}184 & 4.8\% & 4.88 & +1{,}106 & 4.85\\
2{,}000 & 416 & 164 & 9.04 & +63 & +606 & 3.0\% & 4.62 & +586 & 4.61\\
5{,}000 & 1{,}040 & 399 & 9.02 & +28 & +427 & 1.8\% & 4.58 & +440 & 4.57\\
\bottomrule\end{tabular}

\begin{tablenotes}\footnotesize
\item Training solves: dispatch problems solved during training (16 days per step, 13 refits). Time: wall-clock
time on one CPU core, including the level-correction step. Training objective: SPO+ plus squared error
(thousand rupees per day), on each refit's training days, averaged over the 13 refits. Level errors are means of (daily mean forecast before projection) minus
(daily mean actual price), in Rs/MWh. Regret is for the projected forecast, in Rs lakh per MW-year. The
exact level correction is estimated at each refit from training data only. For reference, validation-year
regret is 4.45 for the regression, 4.86 for anchored SPO+ (V3) and 4.08 for EW3.
\end{tablenotes}
\end{threeparttable}
\end{table}

\begin{table}[!htbp]\centering\small
\begin{threeparttable}
\caption{All DFL variants with 150 and 2{,}000 training steps per refit (seed 0). Regret in Rs lakh per
MW-year; EW3 has \BudgetEWval\ (validation) and \BudgetEWhold\ (holdout).}\label{tab:budget}
\footnotesize\setlength{\tabcolsep}{2.0pt}
\begin{tabular}{lccccccccc}\toprule
& \multicolumn{2}{c}{Validation regret} & \multicolumn{2}{c}{Holdout regret} & \makecell{Holdout, 2{,}000 steps,\\minus EW3} & \multicolumn{2}{c}{\makecell{Level error, validation\\(Rs/MWh)}} & \multicolumn{2}{c}{\makecell{Blocks above cap,\\validation (\%)}}\\
\cmidrule(lr){2-3}\cmidrule(lr){4-5}\cmidrule(lr){7-8}\cmidrule(lr){9-10}
Method & 150 & 2{,}000 & 150 & 2{,}000 & [95\% CI] & 150 & 2{,}000 & 150 & 2{,}000\\\midrule
V1 SPO+, unconstrained & 6.60 & 4.62 & 3.10 & 3.11 & $+0.68$ [$+0.54$, $+0.81$] & +2449 & +606 & 8.5 & 3.0\\
V2 SPO+, after projection & 6.95 & 4.93 & 3.15 & 3.16 & $+0.73$ [$+0.59$, $+0.87$] & +2637 & +821 & 10.7 & 5.8\\
V3 SPO+, anchored & 4.90 & 4.86 & 3.10 & 3.09 & $+0.67$ [$+0.53$, $+0.80$] & +344 & +344 & 2.5 & 2.1\\
V4 neural-network SPO+ & 6.43 & 6.97 & 2.89 & 3.10 & $+0.67$ [$+0.49$, $+0.87$] & +344 & +344 & 2.2 & 2.4\\
V5 Fenchel--Young & 5.54 & 5.74 & 3.49 & 3.66 & $+1.23$ [$+1.07$, $+1.39$] & +344 & +344 & 0.9 & 0.4\\
\bottomrule\end{tabular}

\begin{tablenotes}\footnotesize
\item The protocol was fixed before the runs. The 150-step values are the seed-0 runs of the main study; V4 with seeds 1 and 2 is reported in the text; level
errors and blocks above the cap at 150 steps use the seed-0 forecasts of the rerun on the extended data
(Section~\ref{sec:fresh}), which equal the main runs except for the neural network. Level error: mean of
(daily mean forecast before projection minus daily mean price). Holdout interval: circular 7-day block
bootstrap.
\end{tablenotes}
\end{threeparttable}
\end{table}

Table~\ref{tab:decomp} uses the same V1 forecasts (150 steps,
seeds 1 and 2) in different ways in the validation year. Relative to the regression (4.45), V1's extra regret
is 2.07. Removing projection lowers V1's regret from 6.52 to 5.64, a change of 0.89 [0.24, 1.75]; giving the V1
forecasts the regression's daily level lowers it to 4.93, with or without projection. The level and projection
together therefore account for 1.59 of the 2.07 (77\%), and the forecast shape for the remaining 0.48. How the
1.59 is split depends on the order of the changes: removing projection first attributes 0.89 to projection and
0.71 to the level, while correcting the level first leaves nothing for projection. This is a counterfactual
decomposition, not a unique causal split; what is clear is that projection is costly here only because the
level is too high.

\begin{table}[!htbp]\centering\small
\caption{Validation-year regret of the same trained V1 forecasts (150 steps), used in different ways (Rs
lakh per MW-year; seeds 1 and 2).}\label{tab:decomp}
\begin{tabular}{lc}\toprule
How the V1 forecast is used & Regret\\\midrule
Projected to the cap (as tested) & 6.52\\
Not projected & 5.64\\
V1 shape with the regression's daily level, projected & 4.93\\
V1 shape with the regression's daily level, not projected & 4.93\\
Regression (reference) & 4.45\\
\bottomrule\end{tabular}
\end{table}

Table~\ref{tab:sat} groups the same days by the share of V1's
forecast blocks that are above the cap. Projection costs nothing on days with up to 10\% of blocks clipped,
and almost nothing up to 25\%. Above 25\% the cost rises steeply. Theorem~\ref{thm:safe} and Corollary~\ref{cor:marginal} explain this pattern. The battery discharges in about
16 to 23 blocks a day. As long as the block where it runs at partial power is not clipped, clipping its
full-power discharge blocks leaves the raw schedule optimal. When more blocks are clipped than the battery
discharges in, clipping reaches the blocks that matter at the margin.

\begin{table}[!htbp]\centering\small
\begin{threeparttable}
\caption{Cost of projection by the share of forecast blocks above the cap (V1, seeds 1 and 2, all
\IEXnDays{} days).}\label{tab:sat}
\begin{tabular}{lccc}\toprule
Share of blocks above the cap & Days (both seeds) & \makecell{Cost of projection\\(Rs lakh per MW-year)} & \makecell{Cost as \% of\\oracle profit}\\\midrule
None & 1{,}421 & 0.00 & 0.0\\
0 to 10\% & 1{,}182 & 0.00 & 0.0\\
10 to 25\% & 318 & 0.04 & $<0.1$\\
25 to 50\% & 119 & 4.19 & 4.4\\
Over 50\% & 4 & 17.15 & 14.4\\
\bottomrule\end{tabular}
\begin{tablenotes}\footnotesize
\item Cost of projection = regret of the projected V1 forecast minus regret of the same forecast used
without projection, averaged over the days in each group and scaled to a year. It can be negative on
individual days.
\end{tablenotes}
\end{threeparttable}
\end{table}

When the cap fell from Rs\,12{,}000 to Rs\,10{,}000/MWh on 4~April~2023,
prices reached it more often, but the cost of projection did not rise: it fell for V1, whose early level drift
faded, and stayed near zero for the other forecasts (Section~\ref{app:extra}). This descriptive evidence is
consistent with the theory: projection is costly when the forecast level is too high.

For each day and V1 forecast (seeds 1 and 2), we solved the LP at the
raw forecast, took the energy values $\mu_t$ from the solver, and checked the conditions of
Theorem~\ref{thm:safe} (Table~\ref{tab:safe}). The thresholds explain every raw schedule exactly, and no
violation was observed: on all 2{,}312 day-seeds where the conditions hold, the raw schedule is still optimal
after projection, and the two-sided bound of Proposition~\ref{prop:margin}(iii) holds on all 3{,}044
day-seeds. The conditions are conservative. On days with up to 25\% of blocks clipped they fail on about 40\%
of days, mostly (518 days) because of the condition for idle blocks, yet the raw schedule stopped being optimal
on only 2 of 1{,}500 such days. When more than 25\% of blocks are clipped, clipping almost always reaches a
partial-power block (94\%), and the raw schedule stops being optimal on 40\% of days; this is where the cost of
projection in Table~\ref{tab:sat} comes from.

\begin{table}[!htbp]\centering\small
\begin{threeparttable}
\caption{Check of Theorem~\ref{thm:safe} on the V1 forecasts (seeds 1 and 2, all \IEXnDays{} days), by share of
forecast blocks above the cap.}\label{tab:safe}
\setlength{\tabcolsep}{4pt}\begin{tabular}{lccccc}\toprule
Share of blocks & Days & \makecell{Conditions of\\Theorem~\ref{thm:safe} hold (\%)} & \makecell{A partial-power block\\is clipped (\%)} & \makecell{Raw schedule\\not optimal after\\projection (\%)} & \makecell{Discharge blocks\\per day}\\\midrule
None & 1{,}421 & 100.0 & 0.0 & 0.0 & 16.1\\
0 to 10\% & 1{,}182 & 61.2 & 0.0 & 0.0 & 17.5\\
10 to 25\% & 318 & 51.6 & 29.2 & 0.6 & 20.0\\
25 to 50\% & 119 & 3.4 & 94.1 & 39.5 & 20.5\\
Over 50\% & 4 & 0.0 & 100.0 & 100.0 & 23.0\\
\bottomrule\end{tabular}
\begin{tablenotes}\footnotesize
\item The fifth column counts days on which the raw schedule is not a best schedule at the projected forecast. On
the other days it is optimal after projection, although a solver may still return a different tied schedule. Energy
values are the duals returned by the LP solver; other valid duals could make the conditions hold on more
days. ``Days'' counts both seeds.
\end{tablenotes}
\end{threeparttable}
\end{table}

\begin{figure}[!htbp]\centering
\includegraphics[width=0.92\linewidth]{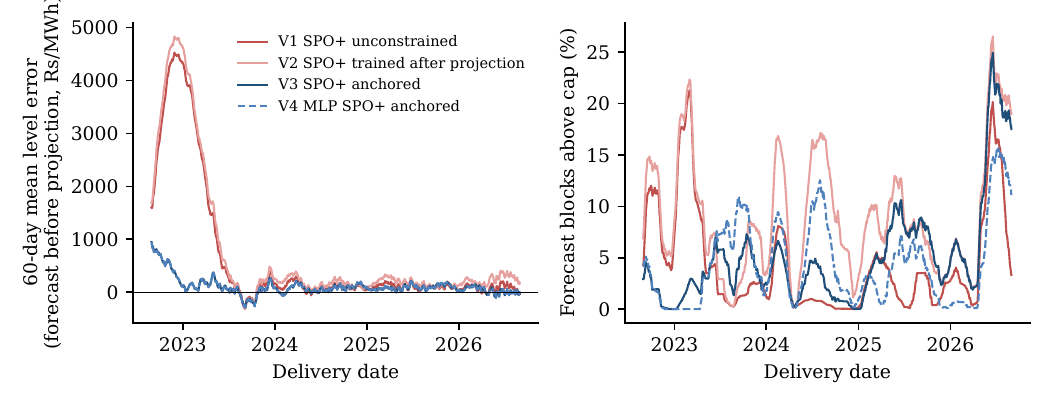}
\caption{India, tuned methods, seed 1. Left: 60-day average level error of the forecast before projection.
Right: share of forecast blocks above the cap. Unconstrained SPO+ and SPO+ trained after projection move
upward early in the sample. The anchored methods keep the regression's level.}\label{fig:drift}
\end{figure}

Training on projected forecasts did not stop the shift: V2 shifts at least as much as V1 (mean level error
+755 against +643 Rs/MWh, and 9.2\% against 4.3\% of forecast blocks above the cap; Figure~\ref{fig:drift}),
consistent with Proposition~\ref{prop:plateau} and Remark~\ref{rem:projpen}. The anchored methods take their
level from the regression, and their results depend little on the learning rate (validation regret 4.90 to
4.98 across all settings, against 6.60 to 32.48 for V1).

In the German test (Section~\ref{app:de}), Assumption~\ref{ass:over} also holds (4.08 against 2.58
MWh per day), and unconstrained SPO+ moved its level up by only about 3 EUR/MWh. No price was cut, and V1 performed the same as its anchored
version.

\subsection{Main comparison}\label{sec:main}
\begin{figure}[!htbp]\centering
\includegraphics[width=0.8\linewidth]{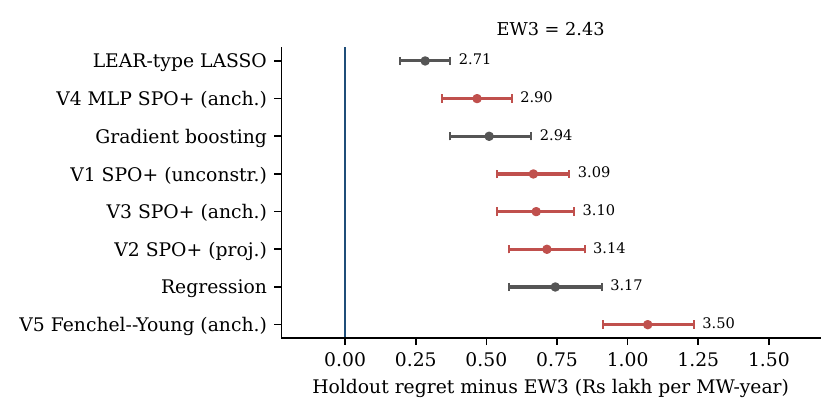}
\caption{Holdout regret minus the regret of EW3, with 95\% bootstrap intervals. Red: DFL methods (average
of three seeds). Grey: statistical forecasts. Labels give the holdout regret (Rs lakh per
MW-year).}\label{fig:holdout}
\end{figure}

Table~\ref{tab:iex} gives the regret of all methods in the validation year and on the holdout. On the holdout, every DFL method is worse than EW3, with every interval excluding zero, by at least \IEXmlpVsEw\ lakh per MW-year. The best
DFL method, the anchored neural-network SPO+ (V4), outperforms the regression it starts from by 0.28
[0.07, 0.48]. It is about as good as gradient boosting and worse than the LASSO. On the holdout, the
linear SPO+ methods are about as good as the regression, and the Fenchel--Young method is the worst.
Different random seeds change the results by at most 0.03. Figure~\ref{fig:holdout} shows the holdout
comparison.

Choosing the DFL settings with hindsight, on the holdout itself, does not change the result. The best anchored linear
SPO+ then reaches 3.09 and the best neural-network SPO+ 2.89, against 2.43 for EW3.

Adding DFL forecasts to the average does not help much:
\begin{itemize}[leftmargin=1.2em]
\item EW3 plus the neural-network SPO+ gives 2.39, a change of $-0.04$ [$-0.09$, $+0.01$];
\item using the neural-network SPO+ instead of the regression gives 2.38, a change of $-0.04$
[$-0.11$, $+0.02$].
\end{itemize}
Weights chosen from past profit or past squared error also do not outperform equal weights. This matches the
known finding that equal weights are hard to beat \citep{smith2009simple}.

\begin{table}[!htbp]\centering\small
\begin{threeparttable}
\caption{India: regret of the 24-hour schedule (Rs lakh per MW-year). Validation: 2 July 2022 to 30 June
2023; holdout: 1 July 2023 to 31 August 2026. DFL rows are the average of three seeds.}\label{tab:iex}
\footnotesize\setlength{\tabcolsep}{3.5pt}
\begin{tabular}{lccccc}\toprule
Method & Validation & Holdout & Seed range & \makecell{Holdout $-$ EW3\\{}[95\% CI]} & \makecell{Holdout $-$ regression\\{}[95\% CI]}\\\midrule
\textbf{EW3} & 4.08 & \textbf{2.43} & -- & -- & $-0.74$ [$-0.91$, $-0.58$]\\
LEAR-type LASSO & 5.24 & 2.71 & -- & $+0.28$ [$+0.19$, $+0.37$] & $-0.46$ [$-0.64$, $-0.28$]\\
Gradient boosting & 5.69 & 2.94 & -- & $+0.51$ [$+0.37$, $+0.66$] & $-0.23$ [$-0.47$, $+0.00$]\\
Regression & 4.45 & 3.17 & -- & $+0.74$ [$+0.58$, $+0.91$] & --\\
Seasonal naive & 5.43 & 4.14 & -- & $+1.71$ & $+0.97$\\\midrule
V1 SPO+, unconstrained & 6.55 & 3.09 & 3.09--3.10 & $+0.67$ [$+0.54$, $+0.79$] & $-0.08$ [$-0.19$, $+0.03$]\\
V2 SPO+, trained after projection & 6.95 & 3.14 & 3.13--3.15 & $+0.71$ [$+0.58$, $+0.85$] & $-0.03$ [$-0.16$, $+0.09$]\\
V3 SPO+, anchored & 4.86 & 3.10 & 3.10--3.11 & $+0.68$ [$+0.54$, $+0.81$] & $-0.07$ [$-0.19$, $+0.06$]\\
V4 neural-network SPO+, anchored & 6.61 & 2.90 & 2.88--2.91 & $+0.47$ [$+0.34$, $+0.59$] & $-0.28$ [$-0.48$, $-0.07$]\\
V5 Fenchel--Young, anchored & 5.56 & 3.50 & 3.49--3.50 & $+1.07$ [$+0.91$, $+1.23$] & $+0.33$ [$+0.19$, $+0.46$]\\
\bottomrule\end{tabular}
\end{threeparttable}
\end{table}

\subsection{Different battery designs}\label{sec:designs}
Table~\ref{tab:robust} changes the battery's duration, efficiency, wear cost and start-and-end level. For
each design we recompute the oracle and all schedules, and we retrain the best DFL method (V4). The
ranking does not change. EW3 has lower holdout regret than V4 in every design, and the interval excludes
zero in \RobustSig\ of 8 designs. The gap is between 0.21 lakh per MW-year (double wear cost) and 0.60
(battery half full at the start and end of the day). EW3 is also better than each of the three forecasts
it averages, in every design.

\begin{table}[!htbp]\centering\small
\begin{threeparttable}
\caption{India, holdout regret for different battery designs (Rs lakh per MW-year).}\label{tab:robust}
\footnotesize\setlength{\tabcolsep}{4pt}
\begin{tabular}{lccccccc}\toprule
Configuration & Oracle & Regression & GB & LASSO & EW3 & V4 MLP SPO+ & EW3 $-$ V4 [95\% CI]\\\midrule
Base (4h, $\eta=0.85$) & 66.7 & 3.17 & 2.94 & 2.71 & \textbf{2.43} & 2.90 & $-0.47$ [$-0.61$, $-0.34$]\\
Duration 2h & 41.6 & 2.53 & 2.47 & 2.14 & \textbf{1.93} & 2.31 & $-0.38$ [$-0.54$, $-0.22$]\\
Duration 8h & 80.7 & 3.32 & 3.07 & 2.95 & \textbf{2.80} & 3.26 & $-0.45$ [$-0.62$, $-0.29$]\\
$\eta=0.80$ & 62.1 & 3.27 & 3.00 & 2.83 & \textbf{2.47} & 2.89 & $-0.43$ [$-0.58$, $-0.29$]\\
$\eta=0.95$ & 75.3 & 2.92 & 2.89 & 2.62 & \textbf{2.34} & 2.77 & $-0.43$ [$-0.55$, $-0.32$]\\
No degradation cost & 86.3 & 2.96 & 2.98 & 2.62 & \textbf{2.37} & 2.91 & $-0.54$ [$-0.67$, $-0.42$]\\
Degradation cost $\times2$ & 50.6 & 4.06 & 3.50 & 3.50 & \textbf{3.14} & 3.34 & $-0.21$ [$-0.40$, $-0.02$]\\
$\bar x=B/2$ & 59.5 & 4.19 & 3.92 & 3.50 & \textbf{3.26} & 3.86 & $-0.60$ [$-0.80$, $-0.38$]\\
\bottomrule\end{tabular}

\begin{tablenotes}\footnotesize
\item V4 is retrained for each design with the settings chosen for the base design (seed 0). The base row therefore differs slightly from the three-seed average in Table~\ref{tab:iex}. ``$\bar x=B/2$'' means the
battery is half full at the start and end of the day.
\end{tablenotes}
\end{threeparttable}
\end{table}

\subsection{A fresh test on new Indian data}\label{sec:fresh}
After the analysis above was complete, we obtained IEX day-ahead prices for 1~September to 24~October 2026.
Before computing any forecast or regret for these days, we fixed a test design (Section~\ref{sec:testdesign}):
the same battery, methods and settings, three seeds and no re-tuning; the same evaluation rule and statistics;
a primary comparison of EW3 with the best DFL method (V4); and a reading rule, under which the ranking is
``confirmed'' if EW3 has lower mean regret than every DFL method, and ``significantly confirmed'' against a
method only if the 95\% interval excludes zero. Because the evaluation rule needs day $D+1$, the test covers
the \FreshN\ delivery days from 1~September to 23~October 2026. Prices for 1--27~September were already in
our original, unused data files; the rest are new. Rerunning the whole pipeline on the extended data
reproduces all earlier forecasts exactly, except the neural network, whose earlier holdout regret changes by
at most 0.005.

This period is very different from the main holdout: 62\% of blocks clear at the cap, against 17\%, and the
oracle profit is Rs\,\FreshOracle\ lakh per MW-year. Table~\ref{tab:fresh} gives the results. EW3 has lower
regret than V4, but the difference, \FreshVfourDiff\ [\FreshVfourCI], is not significant, and although EW3 has
lower mean regret than every DFL method, it is significantly better only than V5. The regression alone has the
lowest regret (\FreshReg, against \FreshEW\ for EW3), although the difference is not significant, and four of
the five DFL methods are significantly worse than the regression they start from (for V2 the interval just
touches zero). So in this fresh, heavily capped period no DFL method outperforms the regression or the
equal-weight average, but the data do not single out EW3 as the best statistical method, and the intervals are
wide.

On the fresh days, with long training windows, V1's level no longer drifts upward (mean
level error \FreshVoneLevel\ Rs/MWh), and its cost of projection is small: \FreshProjCost\ [\FreshProjCostCI].
V2, which is trained after projection, is the only method with an upward level (\FreshVtwoLevel\ Rs/MWh), with
\FreshVtwoAbove\% of its forecast blocks above the cap, consistent with Proposition~\ref{prop:plateau}. No
violation of Theorem~\ref{thm:safe} or of the bound of Proposition~\ref{prop:margin}(iii) was observed. For V1,
the raw schedule stays optimal after projection on every fresh day; the small cost of projection comes from
ties, when the solver picks an equally good schedule that does worse at actual prices. This is the case
covered by the selection clause of Theorem~\ref{thm:safe}.

\begin{table}[!htbp]\centering\small
\begin{threeparttable}
\caption{Fresh test: regret of the 24-hour schedule on \FreshN\ new delivery days, 1 September to 23
October 2026 (Rs lakh per MW-year). The design was fixed before any forecast for these days was computed.
DFL rows are the average of three seeds.}\label{tab:fresh}
\footnotesize\setlength{\tabcolsep}{3.5pt}
\begin{tabular}{lcccc}\toprule
Method & Regret & Seed range & \makecell{Regret $-$ EW3\\{}[95\% CI]} & \makecell{Regret $-$ regression\\{}[95\% CI]}\\\midrule
Regression & 2.67 & -- & $-0.24$ [$-1.28$, $+0.49$] & --\\
\textbf{EW3} & 2.90 & -- & -- & $+0.24$ [$-0.49$, $+1.28$]\\
LEAR-type LASSO & 2.90 & -- & $0.00$ [$-0.65$, $+0.61$] & $+0.23$ [$-0.42$, $+0.96$]\\
Gradient boosting & 3.70 & -- & $+0.80$ [$+0.20$, $+1.38$] & $+1.03$ [$+0.10$, $+2.27$]\\
Seasonal naive & 3.53 & -- & $+0.63$ [$+0.21$, $+1.08$] & $+0.87$ [$+0.09$, $+1.81$]\\\midrule
V3 SPO+, anchored & 3.12 & 3.05--3.22 & $+0.22$ [$-0.39$, $+1.00$] & $+0.45$ [$+0.01$, $+1.00$]\\
V1 SPO+, unconstrained & 3.13 & 2.92--3.45 & $+0.23$ [$-0.41$, $+0.94$] & $+0.47$ [$+0.03$, $+1.03$]\\
V4 neural-network SPO+, anchored & 3.32 & 2.83--4.07 & $+0.42$ [$-0.45$, $+1.35$] & $+0.65$ [$+0.08$, $+1.33$]\\
V2 SPO+, trained after projection & 3.46 & 3.45--3.47 & $+0.56$ [$-0.14$, $+1.45$] & $+0.80$ [$0.00$, $+1.73$]\\
V5 Fenchel--Young, anchored & 4.65 & 4.64--4.66 & $+1.75$ [$+0.23$, $+3.59$] & $+1.99$ [$+0.27$, $+4.00$]\\
\bottomrule\end{tabular}

\begin{tablenotes}\footnotesize
\item Intervals: circular block bootstrap with 7-day blocks and 4{,}000 repetitions over the \FreshN\ days.
\end{tablenotes}
\end{threeparttable}
\end{table}

\subsection{Using the certificate in operation}\label{sec:rule}
What should an operator do with a forecast that will be clipped? The comparisons summarized here are
exploratory: the rules were designed after the forecasting results were known and are evaluated on periods
already used in the study, although each rule was fixed before its own comparison. Section~\ref{app:rules}
gives the rules and the full results (Tables~\ref{tab:policy} and~\ref{tab:select}).

A rule that anchors a forecast only when the certificate of Theorem~\ref{thm:safe} fails lowers the first-year
regret of the unanchored variants V1 and V2 from 6.55 to 5.62 and from 6.95 to 5.72. Anchoring on every day
lowers it further, to 4.95 and 5.28, because a high level costs money even when the raw schedule remains optimal
after clipping. On the holdout the rule changes regret by at most 0.02. For the statistical forecasts the
certificate holds on 92\% to 100\% of days and the rule changes almost nothing, but in the fresh test anchoring
the LASSO to the regression's level raised its regret by 0.67 [0.04, 1.53].

The raw schedule stays optimal after clipping on 97.8\% to 100\% of days for every method. Keeping it whenever
it remains optimal (Corollary~\ref{cor:select}) makes the deployed regret equal to the raw regret on those days.
The effect on average regret is small and goes both ways (between $-0.04$ and $+0.03$ on the holdout), and in the
fresh test it lowers V1's regret by 0.15 [0.00, 0.33]. On the remaining days, choosing the optimal schedule that
the raw forecast prefers lowers first-year regret by 0.40 [0.13, 0.75] for V1 and 0.55 [0.22, 0.95] for V2. Even
with these rules, the best DFL variant stays behind EW3 on the holdout, by 0.46 [0.34, 0.58].

In practice, then, the level of a decision-focused forecast should be anchored, the level of a well-calibrated
statistical forecast should be left alone, the certificate shows on which days clipping may change the schedule,
and ties should be broken in favor of the raw schedule.

\subsection{Size of the gains}\label{sec:gains}
For a 1\,MW/4\,MWh battery in India, EW3 earns Rs\,0.74 lakh per year more than the regression on the
holdout, and Rs\,0.47 lakh more than the best tuned DFL method. These are 23\% of the regression's
holdout regret and 16\% of the best DFL method's. In the first year, unanchored DFL would have lost
Rs\,\IEXvOneLoss\ lakh more than anchored DFL (validation-year regret of \IEXvOneVal\ against
\IEXvThreeVal\ lakh per MW-year). These gains are meaningful relative to the regret a method can reduce
(Proposition~\ref{prop:decomp}), but small relative to the value of the battery: the advantage of EW3 over
the best DFL method is about 0.7\% of the perfect-foresight profit on the holdout (Rs\,66.7 lakh per
MW-year). The value of perfect information is the same for every method.
For a 100\,MW plant, the gain of EW3 over the regression would be about Rs\,0.74 crore per year. This is only a rough figure, because a plant that large would affect prices.

\section{Implications and limitations}\label{sec:impl}
This section brings the theoretical and empirical results together, discusses what they mean for practice,
and states the limitations of the study.

\subsection{Summary of the evidence}
The results separate four effects. Under Assumption~\ref{ass:over}, the losses favor a forecast level that is
too high, but with the squared-error term used in training this effect is small (Theorem~\ref{thm:bias},
Corollary~\ref{cor:reg}, Table~\ref{tab:eta}). The much larger level error early in the sample comes mostly from
limited training: more steps per refit reduce it from about 2{,}450 to about 430 Rs/MWh (Table~\ref{tab:conv}),
and a larger budget for all variants leaves the ranking unchanged (Table~\ref{tab:budget}). A high level hurts
mainly after projection flattens the forecast peaks. When the conditions of Theorem~\ref{thm:safe} hold, the raw
schedule stays optimal, and in our data projection was costly mostly when more than about a quarter of the
blocks were clipped (Tables~\ref{tab:decomp}, \ref{tab:sat} and~\ref{tab:safe}). Out of sample, no DFL variant did
better than the equal-weight statistical average on the holdout or in any of the eight battery designs, and none
did better than the regression or the average in the fresh test (Tables~\ref{tab:iex}, \ref{tab:robust}
and~\ref{tab:fresh}). More training deals with the second effect and anchoring with the third, but neither was
enough for the tested DFL variants to beat the statistical benchmarks. Table~\ref{tab:robsum} lists the
robustness checks and separates the results that follow from the theory from those that depend on the
implementations we tested.

\begin{table}[!htbp]\centering\small
\begin{threeparttable}
\caption{Summary of the robustness checks.}\label{tab:robsum}
\footnotesize\setlength{\tabcolsep}{3pt}
\begin{tabular}{p{3.6cm}p{1.9cm}p{6.1cm}p{2.6cm}}\toprule
Check & Where & Result & Depends on\\\midrule
Random seeds (three per DFL variant) & Table~\ref{tab:iex} & Holdout regret varies by at most 0.03 across seeds & Tested implementations\\
Training budget (2{,}000 steps per refit) & Table~\ref{tab:budget} & Every variant remains worse than EW3 on the holdout, by at least 0.67, all intervals excluding zero; early drift falls sharply & Tested implementations\\
DFL settings chosen with hindsight & Section~\ref{sec:results} & Best DFL variant 2.89 against 2.43 for EW3 & Tested implementations\\
Battery designs (eight) & Table~\ref{tab:robust} & EW3 below the best DFL variant in every design & Tested implementations\\
Second market (Germany) & Table~\ref{tab:de} & Same point ranking; gap to the best DFL variant not significant & Tested implementations\\
Fresh period (\FreshN\ days) & Table~\ref{tab:fresh} & No DFL variant below the regression or EW3 & Tested implementations\\
Direction of the level shift & Tables~\ref{tab:eta}, \ref{tab:de} & Assumption~\ref{ass:over} holds in both markets & Theorem~\ref{thm:bias}: any model trained on these losses\\
Share of clipped blocks & Tables~\ref{tab:sat}, \ref{tab:safe} & Projection costly mainly when more than about 25\% of blocks are clipped & Theorem~\ref{thm:safe}, Corollary~\ref{cor:marginal}\\
Safe-clipping certificate & Table~\ref{tab:safe} & No violation observed in the tested instances & Theorem~\ref{thm:safe}: any model\\
Selection among tied schedules & Table~\ref{tab:select} & Holdout regret changes by $-0.04$ to $+0.03$ & Corollary~\ref{cor:select}: any model\\
Cap change of April 2023 & Section~\ref{app:extra} & Cost of projection follows the level error, not the cap (descriptive) & Theory; descriptive only\\
\bottomrule\end{tabular}
\end{threeparttable}
\end{table}

\subsection{Implications for practice}
Before using a decision-focused forecast, it is worth checking which directions the training loss barely sees.
For a battery this is the price level; for other problems with a conservation constraint it is the conserved
total (Theorem~\ref{thm:general}). Any later step that is not neutral in that direction, such as clipping to
price limits, can turn an error in it into lost profit.

Without anchoring, a decision-focused forecast can end up with a level that is far too high, especially early
on. A high level acts like an extra buying cost and, after clipping, flattens the forecast peaks; in our data
these two effects explain most of the extra regret of unanchored SPO+ (Table~\ref{tab:decomp}). Taking the level
from a separate model removes the problem at little cost, and in the first year it worked better to do this
every day than only when the certificate failed. Changing the level of a well-calibrated statistical forecast,
on the other hand, can make it worse. When forecasts are clipped, it also helps to keep the raw schedule
whenever it remains optimal (Section~\ref{sec:rule}).

Decision-focused methods should also be compared with several good statistical forecasts and their simple
average, using a design fixed in advance. On the holdout, the anchored neural-network SPO+ did 9\% better than
its own regression but 0.47 lakh per MW-year worse than the average, and in the fresh test the regression did
best.

Two further remarks. When all forecasts lead to the same schedule, that schedule is also optimal for their
average, so averaging can change the decision only on days when the forecasts disagree; in our data, this is
where its gain comes from (Section~\ref{app:extra}). And price caps on their own do not make decision-focused
learning fail. The harm comes from clipping forecasts to the cap, which is the operator's choice.

\subsection{Limitations}
Some of the evidence is less independent than we would like. EW3 and the LASSO came out of earlier exploratory
work on the same data, and exploratory SPO+ runs came before the benchmark protocol was fixed
(Section~\ref{sec:testdesign}). The holdout is therefore a fair test of the DFL variants but not of the choice of
EW3. The fresh test does not depend on these choices, but it covers only \FreshN\ days and does not show that EW3
beats the regression. Its protocol was written after those days had passed, and its first 27 days were already
in our original, unused data files. The operating-rule comparisons are exploratory, and German prices never
reach their limits.

On the theory side, the convexity and throughput-matching results hold in general, but the upward shift needs
Assumption~\ref{ass:over}, which holds in our data but need not hold elsewhere. Corollary~\ref{cor:reg} covers
training before projection (V1) only, and we give no theory of how training evolves. Theorem~\ref{thm:safe} is a
sufficient condition, not a necessary one; Proposition~\ref{prop:margin}(iii) bounds the change in regret, but we
have no lower bound that depends only on the share of clipped blocks.

Finally, published DFL implementations
\citep{smets2025forecasting,smets2026proxies,yi2025perturbed,alkhulaifi2026dfl} and larger networks may perform
differently, and we do not replicate them, although Theorem~\ref{thm:bias} applies to any model trained on these
losses without anchoring. We study a small battery that takes prices as given, in one day-ahead market with a
binding cap, and Theorem~\ref{thm:general} is tested only for batteries.

\section{Conclusion}\label{sec:conc}
Decision-focused learning trains forecasts for the decisions they lead to, but forecasts are often
transformed before they reach the decision. We showed three things for problems with a conservation-type constraint. Decision-focused losses identify the forecast level only weakly. When the helper problem inside the loss has excess throughput, as it does in our data, the loss prefers a higher level. A squared-error term keeps this preference small, but limited training can leave the level much higher.
Projecting such a forecast onto a price cap can then erase profitable spreads and increase regret, even as
prediction error falls. We bounded the change in regret caused by post-processing, gave a dual-threshold
certificate under which a raw-optimal schedule remains optimal after clipping, and showed that anchoring
removes the freely learned level shift. On the Indian day-ahead market, an equal-weight average of three
statistical forecasts outperformed every DFL variant we tested on the holdout. On \FreshN\ later days with
most prices at the cap, no variant outperformed the regression or the average, although most differences
from the average were not significant; a German test gave the same point ranking. In capped markets,
decision-focused forecasts should be anchored before use and benchmarked against strong statistical
forecasts. Natural next steps are to extend the theory to batteries that trade in several markets and to
other conservation-constrained problems, and to test the results over a longer period of new data and in
other markets where caps bind.

\section*{Declarations}
\small
\textbf{Data availability.} IEX day-ahead market data are publicly available \citep{iexDAMMarketData}.
The German data are the public benchmark dataset of \citet{lago2021forecasting}. A replication package with the
raw and processed data, all code, cached results, protocol files and checksums is available from the
corresponding author and will be deposited in a public archive.

\textbf{Competing interests.} The authors declare no competing interests.
\normalsize

\bibliographystyle{plainnat}
\bibliography{refs}

\clearpage
\appendix
\renewcommand{\thesection}{EC.\arabic{section}}
\setcounter{table}{0}\renewcommand{\thetable}{EC.\arabic{table}}
\setcounter{figure}{0}\renewcommand{\thefigure}{EC.\arabic{figure}}
\begin{center}{\Large\bfseries Electronic companion}\end{center}

\section{Proofs}\label{app:proofs}
\textbf{Proof of Lemma~\ref{lem:lp}.} Suppose that an LP solution has $c_t>0$ and $d_t>0$ in some block $t$.
Let $\varepsilon=\min\{c_t,d_t/\eta\}>0$ and change only block $t$:
\begin{equation*}
c_t'=c_t-\varepsilon,\qquad d_t'=d_t-\eta\varepsilon .
\end{equation*}
Both powers stay between zero and their limits. Because $\eta=\eta_c\eta_d$, the change in stored energy in
block $t$ is
\begin{equation*}
\Delta\big(-\eta_c\varepsilon+\eta\varepsilon/\eta_d\big)=\Delta\varepsilon\,(-\eta_c+\eta_c)=0,
\end{equation*}
so the storage path and all storage bounds are unchanged. Less charging saves $\Delta\varepsilon(\pi_t+\omega_b)$
and less discharging loses $\Delta\eta\varepsilon(\pi_t-\omega_s)$, so the profit changes by
\begin{equation*}
\Delta\varepsilon\,(\pi_t+\omega_b)-\Delta\eta\varepsilon\,(\pi_t-\omega_s)
=\Delta\varepsilon\,\big[\pi_t(1-\eta)+\omega_b+\eta\,\omega_s\big],
\end{equation*}
which is positive when $\pi_t>\theta$. The original solution was therefore not optimal, and every optimal LP
solution has $c_td_t=0$ when all prices exceed $\theta$. Such a solution is feasible for the problem with
$c_td_t=0$, whose optimal value is at most the LP value, so it is optimal there too. When $\pi_t<\theta$, the
same change in the opposite direction (adding $\varepsilon$ to $c_t$ and $\eta\varepsilon$ to $d_t$, when the
power limits allow it) increases LP profit. When $\eta=1$, the profit change above equals
$\Delta\varepsilon(\omega_b+\omega_s)>0$ for every price. \hfill$\square$

\textbf{Proof of Lemma~\ref{lem:cons}.} Summing the storage dynamics over the day and using $x_H=x_0$ gives
\begin{equation*}
0=x_H-x_0=\Delta\sum_t\big(\eta_cc_t-d_t/\eta_d\big),\qquad\text{so}\qquad
\Delta\sum_td_t=\eta\,\Delta\sum_tc_t=\eta\,T(z).
\end{equation*}
Since $q(z)=d-c$, this gives $\Delta\1^\top q(z)=\eta T(z)-T(z)=-(1-\eta)T(z)$, and therefore
\begin{equation*}
J(\pi+k\1,z)=J(\pi,z)+k\,\Delta\1^\top q(z)=J(\pi,z)-k(1-\eta)T(z).
\end{equation*}
Since $T(z)=\Delta\sum_tc_t$, subtracting $k(1-\eta)T(z)$ is the same as raising the buying cost from
$\omega_b$ to $\omega_b+k(1-\eta)$. For $\alpha>0$,
\begin{align*}
J(\alpha\pi+k\1,z;\omega_b,\omega_s)
&=\alpha\,J(\pi,z;\omega_b/\alpha,\omega_s/\alpha)-k(1-\eta)T(z)\\
&=\alpha\,J\big(\pi,z;(\omega_b+k(1-\eta))/\alpha,\ \omega_s/a\big),
\end{align*}
and multiplication by $\alpha>0$ preserves the maximizers. \hfill$\square$

\textbf{Proof of Proposition~\ref{prop:bound}.} Write $q^\star=q(z^\star)$ and $\hat q=q(\hat z)$, where $\hat z$ is
optimal at $v=\Pi_{\mathcal B}(\hat p)$. First, $\R^{\mathrm{clip}}\ge0$ by optimality of $z^\star$ at $p$.
Second, optimality of $\hat z$ at $v$ gives $J(v,\hat z)-J(v,z^\star)\ge0$. Adding this to
$\R^{\mathrm{clip}}=J(p,z^\star)-J(p,\hat z)$, the cost terms cancel and
\begin{equation*}
\R^{\mathrm{clip}}(p,\hat p)\le\Delta(p-v)^\top(q^\star-\hat q)\le\Delta\sum_t|e_t|\,|q^\star_t-\hat q_t| .
\end{equation*}
For the second bound, write $p-v=-k\1-e^\circ$. By Lemma~\ref{lem:cons},
\begin{equation*}
\Delta(p-v)^\top(q^\star-\hat q)
=-k\,\Delta\1^\top(q^\star-\hat q)-\Delta(e^\circ)^\top(q^\star-\hat q)
=k(1-\eta)\big(T(z^\star)-T(\hat z)\big)-\Delta(e^\circ)^\top(q^\star-\hat q),
\end{equation*}
and bounding each term by its absolute value gives the second bound. \hfill$\square$

\textbf{Proof of Proposition~\ref{prop:decomp}.} Adding and subtracting $\Pi^\star(\nu)$ and $J(\nu,\hat z)$,
\begin{equation*}
\Pi^\star(p)-J(p,\hat z)=\big[\Pi^\star(p)-\Pi^\star(\nu)\big]+\big[\Pi^\star(\nu)-J(\nu,\hat z)\big]
+\big[J(\nu,\hat z)-J(p,\hat z)\big].
\end{equation*}
Since $J(p,\hat z)$ is affine in $p$ and $\hat z$ is $\mathcal I_D$-measurable,
$\E[J(p,\hat z)\mid\mathcal I_D]=J(\nu,\hat z)$, so the last term has zero expectation, and taking expectations
gives the identity. The first term is nonnegative in expectation by Jensen's inequality,
$\E[\Pi^\star(p)\mid\mathcal I_D]\ge\Pi^\star(\nu)$, because $\Pi^\star$ is convex as a maximum of affine
functions. The second term is nonnegative by definition of $\Pi^\star$. \hfill$\square$

\textbf{Properties of SPO+.} Convexity holds because $\Pi^\star$ is convex and $J(\hat p,z^\star)$ is affine in
$\hat p$. At $\hat p=p$,
\begin{equation*}
\ell(p,p)=\Pi^\star(p)-2\Pi^\star(p)+\Pi^\star(p)=0 .
\end{equation*}
For the bound, let $\hat z^{\,\mathrm{raw}}$ be optimal at $\hat p$ itself. The cost term enters
$2J(\hat p,z)-J(p,z)$ with coefficient $2-1=1$, so $2J(\hat p,z)-J(p,z)=J(2\hat p-p,z)$ for every $z$. Using
first $J(\hat p,\hat z^{\,\mathrm{raw}})\ge J(\hat p,z^\star)$ and then the definition of $\Pi^\star$,
\begin{align*}
-J(p,\hat z^{\,\mathrm{raw}})
&\le2J(\hat p,\hat z^{\,\mathrm{raw}})-J(p,\hat z^{\,\mathrm{raw}})-2J(\hat p,z^\star)\\
&=J(2\hat p-p,\hat z^{\,\mathrm{raw}})-2J(\hat p,z^\star)\ \le\ \Pi^\star(2\hat p-p)-2J(\hat p,z^\star).
\end{align*}
Adding $\Pi^\star(p)=J(p,z^\star)$ to both sides gives $\ell(\hat p,p)\ge\R^{\mathrm{raw}}(p,\hat p)$. The argument
needs $\hat z^{\,\mathrm{raw}}$ to be optimal for the \emph{same} vector that enters the loss; it uses only this
optimality and the fact that $J$ is affine in prices, so it does not need $\F$ to be convex. Applied to
$\Pi_{\mathcal B}(\hat p)$, it gives
\begin{equation*}
\ell\big(\Pi_{\mathcal B}(\hat p),p\big)\ge\R^{\mathrm{raw}}\big(p,\Pi_{\mathcal B}(\hat p)\big)=\R^{\mathrm{clip}}(p,\hat p).
\end{equation*}
The subgradient follows from Danskin's theorem. \hfill$\square$

\textbf{Proof of Theorem~\ref{thm:eta1}.} (i) Note that $0\le(1-\eta)T(z)\le(1-\eta)\bar T$. By Lemma~\ref{lem:cons},
\begin{equation*}
\Pi^\star(2\hat p-p+2k\1)=\max_z\{J(2\hat p-p,z)-2k(1-\eta)T(z)\},\qquad
J(\hat p+k\1,z^\star)=J(\hat p,z^\star)-k(1-\eta)T(z^\star).
\end{equation*}
Substituting both into the definition of $\ell$ gives
\begin{equation*}
\ell(\hat p+k\1,p)-\ell(\hat p,p)=\Big[\max_z\{J(2\hat p-p,z)-2k(1-\eta)T(z)\}-\Pi^\star(2\hat p-p)\Big]+2k(1-\eta)T(z^\star).
\end{equation*}
For $k\ge0$, every objective in the maximum falls by between $0$ and $2k(1-\eta)\bar T$, so the bracket lies
in $[-2k(1-\eta)\bar T,0]$, and the last term lies in $[0,2k(1-\eta)\bar T]$. Their sum therefore lies
in $[-2k(1-\eta)\bar T,\,2k(1-\eta)\bar T]$. The case $k<0$ is symmetric.

(ii) If $k\ge\bar p-\min_tp_t$, every block of $p+k\1$ is at least $\bar p$, so the projection is $\bar p\1$. By
Lemma~\ref{lem:cons},
\begin{equation*}
J(\bar p\1,z)=\bar p\,\Delta\1^\top q(z)-C(z)=-\bar p(1-\eta)T(z)-C(z)\le-C(z),
\end{equation*}
and $C(z)>0$ for every schedule that trades. Idling ($c=d=0$, $x_t\equiv\bar x$) is feasible with value 0, so
it is the unique maximizer. It earns 0 at realized prices, so the regret is $\Pi^\star(p)$. The bound on $\ell$
follows from (i) with $\hat p=p$ and $\ell(p,p)=0$.

(iii) A level shift of the model shifts every sample forecast by the same amount. By (i), the empirical risk is
unchanged when $\eta=1$, so every shift of a minimizer is again a minimizer. For
$k\ge\max_i\big(\bar p_i-\min_t\hat p^{(i)}_t\big)$, where $\hat p^{(i)}$ is the forecast for sample day $i$, which is finite for a finite sample, the argument of (ii)
applies on every sample day. \hfill$\square$

\textbf{Proof of Lemma~\ref{lem:mono}.} By Lemma~\ref{lem:cons}, $z_i$ maximizes $J(\pi,z)-k_i(1-\eta)T(z)$, so
\begin{align*}
J(\pi,z_1)-k_1(1-\eta)T(z_1)&\ge J(\pi,z_2)-k_1(1-\eta)T(z_2),\\
J(\pi,z_2)-k_2(1-\eta)T(z_2)&\ge J(\pi,z_1)-k_2(1-\eta)T(z_1).
\end{align*}
Adding the two inequalities, the terms $J(\pi,z_1)$ and $J(\pi,z_2)$ cancel, and
\begin{equation*}
(k_2-k_1)(1-\eta)\big(T(z_1)-T(z_2)\big)\ge0 .
\end{equation*}
Since $k_2>k_1$ and $\eta<1$, this gives $T(z_1)\ge T(z_2)$. \hfill$\square$

\textbf{Proof of Theorem~\ref{thm:bias}.} (i) By Lemma~\ref{lem:cons}, for every realization,
\begin{equation*}
\ell(\hat p+k\1,p)=\max_{z\in\mathcal V}\big\{J(2\hat p-p,z)-2k(1-\eta)T(z)\big\}+2k(1-\eta)T(z^\star)
+\Pi^\star(p)-2J(\hat p,z^\star).
\end{equation*}
The maximum is a maximum of finitely many affine functions of $k$, one per point of $\mathcal V$, so it is
convex in $k$. Its right derivative is the largest slope among the active functions, and its left derivative
the smallest:
\begin{equation*}
-2(1-\eta)\,T^-(2\hat p-p+2k\1)\qquad\text{and}\qquad-2(1-\eta)\,T^+(2\hat p-p+2k\1).
\end{equation*}
The remaining terms are affine in $k$, with slope $2(1-\eta)T(z^\star)$. For convex functions the difference
quotients are monotone. Throughput is bounded by $\bar T$, and by the integrability assumption the
quotients are bounded by an integrable function. So expectation and one-sided differentiation can be exchanged
by monotone convergence, which gives the stated derivatives of $\Phi$.

(ii) By convexity, $k\in K^\star$ if and only if $\Phi'_-(k)\le0\le\Phi'_+(k)$, which is the stated condition.
Monotonicity follows from Lemma~\ref{lem:mono}, applied to each realization. For the existence claim, write
$\pi=2\hat p-p$. As $k\to+\infty$: idling is feasible with value 0, so every best schedule $z$ at $\pi+2k\1$
satisfies
\begin{equation*}
J(\pi,z)-2k(1-\eta)T(z)\ge0,\qquad\text{so}\qquad T(z)\le\frac{\Pi^\star(\pi)}{2k(1-\eta)}\to0 .
\end{equation*}
So, for each realization, $T^+(\pi+2k\1)\to0$ pointwise, and $0\le T^+(\pi+2k\1)\le\bar T$ for every $k$.
The constant $\bar T$ is integrable, so by dominated convergence $\E\,T^+(\pi+2k\1)\to0<\E\,T(z^\star)$, and
$\Phi'_-(k)>0$ for large $k$. As $k\to-\infty$: let $\bar z$ be a schedule with $T(\bar z)=\bar T$. Every
best schedule $z$ satisfies
\begin{equation*}
J(\pi,z)-2k(1-\eta)T(z)\ge J(\pi,\bar z)-2k(1-\eta)\bar T,\qquad\text{so}\qquad
\bar T-T(z)\le\frac{\Pi^\star(\pi)-J(\pi,\bar z)}{2|k|(1-\eta)}\to0 .
\end{equation*}
So $T^-(\pi+2k\1)\to\bar T$ pointwise, with the same bound $0\le T^-\le\bar T$, and by dominated
convergence $\E\,T^-(\pi+2k\1)\to\bar T>\E\,T(z^\star)$; hence $\Phi'_+(k)<0$ for very negative $k$. A
convex function whose right derivative is negative at one point and whose left derivative is positive at a
larger point has a nonempty, closed and bounded set of minimizers.

(iii) Under Assumption~\ref{ass:over}, $\Phi'_+(0)=-2(1-\eta)\big(\E\,T^-(2\hat p-p)-\E\,T(z^\star)\big)<0$. For
$k<0$, convexity gives $\Phi(k)\ge\Phi(0)+\Phi'_+(0)\,k>\Phi(0)$, and $\Phi'_+(0)<0$ means that $\Phi$ decreases
just to the right of 0. So no $k\le0$ minimizes $\Phi$.

(iv) By Lemma~\ref{lem:cons}, the helper problem at $2\hat p-p+2k\1$ has the same best schedules as the problem
at $2\hat p-p$ with the buying cost raised by $\delta=2k(1-\eta)$. So
\begin{equation*}
T^\pm(2\hat p-p+2k\1)=T^\pm_{2k(1-\eta)}(2\hat p-p),
\end{equation*}
and the condition in (ii) holds at $k$ if and only if $2k(1-\eta)\in D^\star$. So
$K^\star=\{\delta/(2(1-\eta)):\delta\in D^\star\}$ and
\begin{equation*}
\min K^\star=\frac{\min D^\star}{2(1-\eta)}\ge\frac{\delta_0}{2(1-\eta)}\to\infty\qquad(\eta\to1).
\end{equation*}

(v) By Lemma~\ref{lem:cons},
\begin{equation*}
L_\sigma(\hat p+k\1,p)=\E_Z\max_z\big\{J(\hat p+\sigma Z,z)-k(1-\eta)T(z)\big\}-J(\hat p,z^\star)+k(1-\eta)T(z^\star),
\end{equation*}
and the same arguments apply with slope $1-\eta$ in place of $2(1-\eta)$. \hfill$\square$

\textbf{Proof of Corollary~\ref{cor:reg}.} With $\bar e=\E[H^{-1}\1^\top(\hat p-p)]$, the penalty is
\begin{equation*}
\lambda\,\E\big[H^{-1}\|\hat p+k\1-p\|^2\big]=\lambda\big(k^2+2k\bar e+\E[H^{-1}\|\hat p-p\|^2]\big),
\end{equation*}
which is strictly convex in $k$ with derivative $2\lambda(k+\bar e)$. Adding the convex $\Phi$ gives a strictly
convex function. It tends to infinity as $|k|\to\infty$, because $\Phi$ is bounded below by 0, so it has a unique
minimizer $k_\lambda$. (i) The minimizer satisfies
\begin{equation*}
\Phi'_-(k_\lambda)+2\lambda(k_\lambda+\bar e)\ \le\ 0\ \le\ \Phi'_+(k_\lambda)+2\lambda(k_\lambda+\bar e).
\end{equation*}
Substituting the derivatives from Theorem~\ref{thm:bias}(i) and dividing by 2 gives the stated condition. (ii)
Both $\E\,T^\pm$ and $\E\,T(z^\star)$ lie in $[0,\bar T]$, so their difference is at most $\bar T$ in
absolute value, and (i) gives the bound. (iii) At $k=-\bar e$ the penalty has zero derivative, so
$\mathcal L'_{\lambda,+}(-\bar e)=\Phi'_+(-\bar e)<0$ by the assumption, and $k_\lambda>-\bar e$. If
$m=\min K^\star>-\bar e$, then
\begin{equation*}
\mathcal L'_{\lambda,+}(m)=\Phi'_+(m)+2\lambda(m+\bar e)>0,
\end{equation*}
because $\Phi'_+(m)\ge0$ at a minimizer of $\Phi$. So $k_\lambda\le m$. \hfill$\square$

\textbf{Proof of Proposition~\ref{prop:margin}.} Since $J(v,z)-J(u,z)=\Delta(v-u)^\top q(z)$ and $v-u$ is zero
outside $\mathcal C$,
\begin{equation*}
J(v,z_u)-J(v,z)=J(u,z_u)-J(u,z)+\Delta\sum_{t\in\mathcal C}(v_t-u_t)\big(q_t(z_u)-q_t(z)\big).
\end{equation*}
(i) Because $J(v,\cdot)$ is linear in the schedule, it attains its maximum over $\F$ at a point of $\mathcal V$, so
$z_u$ is a best schedule at $v$ if and only if the right-hand side is nonnegative for every $z\in\mathcal V$, which
is the stated condition; the set of such $v$ is the normal cone $N(z_u)$, an intersection of finitely many
half-spaces. (ii) The last sum is at least $-\Delta\sum_{t\in\mathcal C}|u_t-v_t|\,|q_t(z_u)-q_t(z)|$, so under the
condition the right-hand side is positive for every other $z\in\mathcal V$. So $z_u$ is the unique best schedule at
$v$, and the two regrets are equal. (iii) Since $\R^{\mathrm{clip}}-\R^{\mathrm{raw}}=J(p,z_u)-J(p,z_v)$ and
$J(p,z)=J(v,z)+\Delta(p-v)^\top q(z)$,
\begin{equation*}
\R^{\mathrm{clip}}-\R^{\mathrm{raw}}=\big[J(v,z_u)-J(v,z_v)\big]+\Delta(p-v)^\top\big(q(z_u)-q(z_v)\big)
\le\Delta(p-v)^\top\big(q(z_u)-q(z_v)\big),
\end{equation*}
because $z_v$ is optimal at $v$. In the same way, using $J(p,z)=J(u,z)+\Delta(p-u)^\top q(z)$ and the optimality of
$z_u$ at $u$,
\begin{equation*}
\R^{\mathrm{clip}}-\R^{\mathrm{raw}}=\big[J(u,z_u)-J(u,z_v)\big]+\Delta(p-u)^\top\big(q(z_u)-q(z_v)\big)
\ge\Delta(p-u)^\top\big(q(z_u)-q(z_v)\big).
\end{equation*}
The difference between the two bounds is
\begin{equation*}
\Delta(u-v)^\top\big(q(z_u)-q(z_v)\big)=\big[J(u,z_u)-J(u,z_v)\big]+\big[J(v,z_v)-J(v,z_u)\big]\ge0,
\end{equation*}
and $u-v$ is zero outside $\mathcal C$. \hfill$\square$

\textbf{Proof of Theorem~\ref{thm:safe}.} Because all prices exceed $\theta$, the LP is exact at $u$ and at $v$
(Lemma~\ref{lem:lp}), and the optimality (KKT) conditions of the LP are necessary and sufficient. Let $\mu$ be the
optimal dual solution for which the conditions of the theorem hold, with $\mu_t$ the multiplier of the storage
equation of block $t$. Since $z_u$ and $\mu$ are optimal primal and dual solutions, they satisfy complementary
slackness together. At a price vector $\pi$, the optimality conditions are of two kinds. For every block $t$,
\begin{equation*}
c_t=\bar c\ \text{if}\ \pi_t<\theta^c_t,\quad c_t=0\ \text{if}\ \pi_t>\theta^c_t,\quad
d_t=\bar d\ \text{if}\ \pi_t>\theta^d_t,\quad d_t=0\ \text{if}\ \pi_t<\theta^d_t,
\end{equation*}
with fractional power only where $\pi_t$ equals the corresponding threshold. For every storage level $x_t$,
\begin{equation*}
\mu_{t+1}-\mu_t\ge0\ \text{if}\ x_t=B,\qquad\mu_{t+1}-\mu_t\le0\ \text{if}\ x_t=0,\qquad
\mu_{t+1}-\mu_t=0\ \text{otherwise}.
\end{equation*}
The storage conditions do not involve prices, and at $v$ prices change only in $\mathcal C$. In a block cut to
the cap, the price falls from $u_t$ to $\bar p$. Full charging remains consistent, because a lower price makes
charging more attractive. Full discharging remains consistent if $\bar p\ge\theta^d_t$. Idling remains consistent
if $\bar p\ge\theta^c_t$, since the price stays below $\theta^d_t$. Fractional power is excluded by (a). The case
of a block raised to the floor is symmetric. So $(z_u,\mu)$ satisfies the optimality conditions at $v$, and $z_u$
is a best schedule at $v$. \hfill$\square$

\textbf{Proof of Corollary~\ref{cor:marginal}.} Fractional discharge at $m$ means $u_m=\theta^d_m$. If the storage
is strictly between its limits from $m$ to $t$, the storage conditions give $\mu$ constant on these blocks, so
\begin{equation*}
\theta^d_t=\mu_t/\eta_d+\omega_s=\mu_m/\eta_d+\omega_s=\theta^d_m=u_m\le\bar p,
\end{equation*}
where the last inequality holds because block $m$ is not clipped. \hfill$\square$

\textbf{Proof of Theorem~\ref{thm:general}.} By Assumption~\ref{ass:cons}, for every $z$,
\begin{equation*}
J(u+kv,z)=J(u,z)+k\,\Delta r^\top q(z)=J(u,z)-k\,aT(z)+kb .
\end{equation*}
In the SPO+ loss this gives
\begin{align*}
\Pi^\star(2u-p+2kv)&=\max_z\big\{J(2u-p,z)-2k\,aT(z)\big\}+2kb,\\
-2J(u+kv,z^\star)&=-2J(u,z^\star)+2k\,aT(z^\star)-2kb,
\end{align*}
so the terms in $b$ cancel. The proofs of Theorem~\ref{thm:eta1}(i), Lemma~\ref{lem:mono},
Theorem~\ref{thm:bias}(i), (iii)--(v) and Corollary~\ref{cor:reg} then apply verbatim with $(1-\eta)T$ replaced
by $aT$. For the existence claim in (ii), write $\pi=2\hat u-p$. Because
$\mathcal V_{\mathcal Z}$ is finite, for $k$ large enough every best decision at $\pi+2kv$ minimizes $T$ over
$\mathcal V_{\mathcal Z}$ (the term $-2kaT$ dominates the bounded differences in $J(\pi,\cdot)$), and for $k$ small
enough every best decision maximizes $T$. So, for each realization, the helper $T$ converges pointwise to
$\min_{\mathcal Z}T$ and to $\bar T$, and $0\le T\le\bar T$ throughout. By dominated convergence, with the
integrable constant bound $\bar T$, the expected helper $T$ tends to $\min_{\mathcal Z}T$ and to $\bar T$. Under
the stated condition the one-sided derivatives therefore change sign, and a convex function with this property
has a nonempty, closed and bounded set of minimizers. For (iii), the derivative of the penalty along $r$ is
\begin{equation*}
\frac{d}{dk}\,\lambda\,\E\big[H^{-1}\|\hat u+kv-p\|^2\big]=\frac{2\lambda}{H}\big(k\|r\|^2+r^\top\E(\hat u-p)\big),
\end{equation*}
so the penalized objective is strictly convex. At its minimizer the one-sided derivatives of $\Phi$ are at most
$2a\bar T$ in absolute value, so
\begin{equation*}
\Big|k_\lambda+\frac{r^\top\E(\hat u-p)}{\|r\|^2}\Big|\le\frac{a\bar TH}{\lambda\|r\|^2}.
\end{equation*}
If $a=0$, the profit of every decision changes by the same constant $kb$, so raw decisions are unchanged. For
(iv), $\ell(p+kv,p)=\ell(p,p)=0$ by (i) with $a=0$, and $g(p+kv)=s$ fixes the deployed forecast. For (v), each
coordinate $g(u)_t$ with $t\in S$ is constant on the stated interval. The other coordinates of $g(u)$ do not
depend on $(u_t)_{t\in S}$, because $g$ acts coordinate by coordinate. So $g(u)$, and with it the loss, is
unchanged, and any chain rule through $g$ multiplies by a zero derivative in these coordinates. For (vi),
\begin{equation*}
(u+kv)-\frac{w^\top(u+kv)}{w^\top r}\,r=u-\frac{w^\top u}{w^\top r}\,r\qquad\text{for every }k,
\end{equation*}
so the anchored forecast does not change when $u$ moves along $r$. \hfill$\square$

\textbf{Proof of Corollary~\ref{cor:select}.} If $z_u$ is a best schedule at $v$, the rule deploys $z_u$, which is
also the schedule deployed without projection, so the two regrets are equal; Theorem~\ref{thm:safe} gives
conditions under which $z_u$ is a best schedule at $v$. Otherwise the rule deploys a best schedule $z_v$ at $v$,
and Proposition~\ref{prop:margin}(iii) applies to this pair $(z_u,z_v)$. \hfill$\square$

\textbf{Proof of Proposition~\ref{prop:plateau}.} If $u_t>\bar p$, then $\Pi_{\mathcal B}(u)_t=\bar p$ for all values
of $u_t$ near the current one, and similarly below the floor. So on $\mathcal C(u)$ the projection is locally
constant in $u_t$, and $\ell(\Pi_{\mathcal B}(u),p)$ does not change with small changes in $(u_t)_{t\in\mathcal
C(u)}$. If $u_t\ge\bar p$ for all $t$, then for all $k\ge0$
\begin{equation*}
\Pi_{\mathcal B}(u+k\1)=\bar p\1,\qquad\text{so}\qquad\ell\big(\Pi_{\mathcal B}(u+k\1),p\big)=\ell(\bar p\1,p).
\end{equation*}
\hfill$\square$

\textbf{Proof of Corollary~\ref{cor:anchor}.} (i) For every $k$,
\begin{equation*}
(u+k\1)-H^{-1}\1^\top(u+k\1)\,\1+\hat L\1=u-H^{-1}(\1^\top u)\,\1+\hat L\1=\tilde u,
\end{equation*}
so $\tilde u$, and with it the training loss, does not depend on a common shift of $u$. (ii) The deployed
forecast $\hat p=\Pi_{\mathcal B}(\tilde u)$ is already projected, so
\begin{equation*}
\ell(\hat p,p)\ge\R^{\mathrm{raw}}(p,\hat p)=\R^{\mathrm{clip}}(p,\tilde u).
\end{equation*}
(iii) Writing $u^\circ$ and $p^\circ$ for the deviations from the daily means,
\begin{equation*}
\tilde u-p=\big(\hat L-H^{-1}\1^\top p\big)\1+\big(u^\circ-p^\circ\big),
\end{equation*}
so the level component of the error of $\tilde u$ is that of $\hat L$, and Proposition~\ref{prop:bound}
applies. \hfill$\square$

\section{Implementation details}\label{app:impl}
\textbf{Dispatch optimization.} Each dispatch problem is first solved as a linear program with HiGHS; if the
solution charges and discharges in the same block, it is re-solved as a mixed-integer program (see below).
During training, the linear-programming model is kept in memory and only the prices change between solves.
For evaluation, each problem is solved from scratch, so that ties are broken in the same way for every
method. Both ways give the same profit to within
$4\times10^{-11}$.

\textbf{No simultaneous charging and discharging.} By Lemma~\ref{lem:lp}, the LP is exact when every price
exceeds $\theta=-(\omega_b+\eta\omega_s)/(1-\eta)$. We solve the LP and, if a solution charges and discharges in
the same block, re-solve with binary variables $y_t$ ($c_t\le\bar c\,y_t$, $d_t\le\bar d(1-y_t)$). For the
Indian base battery ($\theta\approx-6{,}780$ Rs/MWh), all prices exceeded $\theta$ in every audited solve, so the
LP was exact and no re-solve was needed: the oracle and all forecast schedules on all evaluation days; the SPO+
helper problems for the regression and V1--V4; 6{,}088 Fenchel--Young draws; 211{,}200 training-time solves of
V1, V2 and V5 in the validation year (lowest price $-4{,}101$ Rs/MWh); and the level searches of
Table~\ref{tab:eta} (lowest price $-3{,}824$ Rs/MWh, at least 1{,}107 Rs/MWh above $\theta$ at every
efficiency). Realized prices and projected forecasts are never negative, and in the level-correction
bisection of Table~\ref{tab:conv} the final brackets lie where all helper prices exceed $\theta$. With
$\eta=0.80$ ($\theta\approx-4{,}930$), none of 70{,}400 training-time solves was affected. With no wear cost,
$\theta$ is only about $-250$ Rs/MWh, and 0.34\% of training-time helper solves charged and discharged at once;
we reran this design with the exact re-solve in every solve, and 1{,}524 of 145{,}249 solves needed it
(Table~\ref{tab:robust}). In Germany ($\theta\approx-77$ EUR/MWh; 16 evaluation hours below it), the LP oracle
charged and discharged at once on 3 of 2{,}184 days and the gradient boosting schedule on 2 days, so all German
results use the exact re-solve everywhere, including all training-time helper problems; during training it was
needed in 0.013\% to 0.021\% of solves.

\textbf{Code checks.} On 200 randomly chosen Indian days, $\ell(p,p)=0$ exactly, and
$\ell(\hat p,p)-\R^{\mathrm{raw}}(p,\hat p)\ge35.65$ rupees for the projected forecast $\hat p$, so the loss is
always above the regret. With the oracle schedule fixed, the SPO+ subgradient passes the convexity check
(directional derivative at least $g^\top v$) in all 400 random directions, to within $2\times10^{-11}$, and the
derivative with respect to the level matches Theorem~\ref{thm:bias}(i) to machine precision.

\textbf{DFL training.} Model outputs are in units of Rs\,1{,}000/MWh (India) or EUR\,10/MWh (Germany).
The loss is SPO+ (or Fenchel--Young) plus $\lambda$ times the mean squared error of the forecast that is
used. Gradients are set to zero in blocks cut by the projection. For anchored methods, they are also
centered within the day. Adam uses $\beta=(0.9,0.999)$. The first refit takes 400 steps from the starting
model. Each later 28-day refit takes 150 steps. Each step uses 16 training days drawn at random from all
past data. The neural network (MLP) has 20 (India) or 22 (Germany) inputs per block: lagged prices,
summary values of day $D-1$, block, day of week and month terms, and in Germany the load and
wind-plus-solar forecasts. It has one hidden layer of 32 tanh units and is first trained for 3{,}000 steps
on squared error. The Fenchel--Young gradient averages 4 random normal draws per day.

\textbf{Chosen settings.} India: V1 and V2 $(\lambda,\text{lr})=(1,3\times10^{-3})$; V3
$(0.1,10^{-3})$; V4 $(1,10^{-3})$; V5 $(\lambda,\sigma,\text{lr})=(0.1,1000,3\times10^{-4})$.
Germany: \DEsel.

\textbf{Code and data.} The IEX data are public market data, and the German data are the public benchmark
dataset of \citet{lago2021forecasting}. All code (data preparation, dispatch optimization, forecasts, DFL training and evaluation) is included in the replication package that accompanies the submission.

\section{Second test: German day-ahead market}\label{app:de}
\textbf{Purpose of the test.} The main results use one market, and EW3 was first found in exploratory work on
the same data. To check whether the ranking holds elsewhere, we repeat the comparison on a second
market. We set up this test after seeing the Indian results but before seeing any German result. German
prices never reach their limits, so this test checks the ranking and whether the level shift appears. It
cannot check the loss caused by a binding cap.

\textbf{Data and battery.} We use the public benchmark dataset of \citet{lago2021forecasting} for the EPEX
day-ahead market (Germany--Austria). It has hourly prices from 9~January~2012 to 31~December~2017, and
the load forecast and wind-plus-solar forecast published before the auction. The price limits are
$-500$ and $+3{,}000$~EUR/MWh, and they never bind: actual prices range from \DEpmin\ to
\DEpmax~EUR/MWh. A day has $H=24$ blocks with $\Delta=1$. The battery is the same as in India. Its wear
cost is EUR\,12 per MWh cycled (the Indian value in euros) and the fee is EUR\,0.25/MWh for buying and for
selling. With perfect foresight, the battery would earn EUR\,\DEoracle\ per MW-year.

\textbf{Methods and design.} The methods are those of Section~\ref{sec:design}, with three changes:
\begin{itemize}[leftmargin=1.2em]
\item every model also uses the load and wind-plus-solar forecasts for day $D$;
\item the regression also has a weekend dummy and day-of-year terms;
\item the LASSO also uses these forecasts for days $D$, $D-1$ and $D-7$.
\end{itemize}
Because the limits never bind, V2 is the same as V1. We therefore run V1, V3 and V4. We drop V5 because it was
the weakest method in India. The validation year is 2015 and the holdout is 2016--2017 (731 days), the
same test period as in \citet{lago2021forecasting}. Tuning, seeds and statistics are as in India.

\begin{figure}[!htbp]\centering
\includegraphics[width=0.75\linewidth]{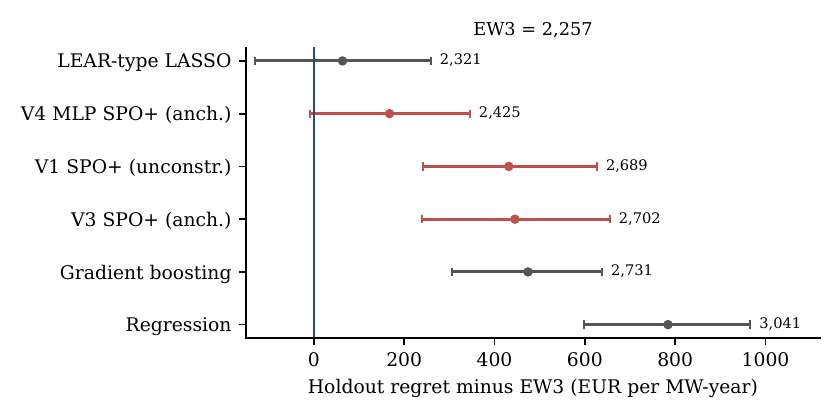}
\caption{Germany: holdout regret minus the regret of EW3, with 95\% bootstrap intervals. Red: DFL methods
(average of three seeds). Grey: statistical forecasts. Labels give the holdout regret (EUR per
MW-year).}\label{fig:holdout_de}
\end{figure}

\textbf{Results.} Table~\ref{tab:de} and Figure~\ref{fig:holdout_de} give the results. The ranking is the
same as in India.
\begin{itemize}[leftmargin=1.2em]
\item EW3 has the lowest holdout regret, EUR\,\DEewH\ per MW-year. The LASSO is not significantly
different from it.
\item The best DFL method is again the anchored neural-network SPO+, at \DEmlpH. This is \DEbestVsEw\
[\DEbestVsEwCI] above EW3. The point estimate favors EW3, but the interval includes zero, so in Germany
this gap is not statistically significant at the 5\% level. The linear SPO+ methods are about 430 to 450
above EW3, and these gaps are significant.
\item All three DFL methods outperform the regression they start from, by 339 to 617 EUR per MW-year. (In
India, only the neural-network method does.) So, without a binding cap, DFL improves on its starting model, but it does not outperform an average of good statistical forecasts.
\end{itemize}

\textbf{Level shift.} At the regression forecast, the SPO+ helper problem
charges 4.08 MWh per day, while the oracle charges 2.58. So Theorem~\ref{thm:bias} predicts an upward
shift. Unconstrained SPO+ did move its level up compared with its starting regression, by \DEshiftVal\
EUR/MWh in the validation year. But the shift was small compared with German prices, no hour ever went
outside the price limits, and V1 and its anchored version perform the same on the holdout (\DEvOneH\ and
\DEvThreeH\ EUR per MW-year). The upward shift is observed, but its effect on regret is not statistically significant here.

\begin{table}[!htbp]\centering\small
\begin{threeparttable}
\caption{Germany: regret of the 24-hour schedule (EUR per MW-year). Validation: 2015; holdout: 2016--2017.
DFL rows are the average of three seeds.}\label{tab:de}
\footnotesize\setlength{\tabcolsep}{3pt}
\begin{tabular}{lcccccc}\toprule
Method & \makecell{RMSE\\(holdout)} & Validation & Holdout & Seed range & \makecell{Holdout $-$ EW3\\{}[95\% CI]} & \makecell{Holdout $-$ regression\\{}[95\% CI]}\\\midrule
\textbf{EW3} & 6.71 & 2{,}024 & \textbf{2{,}257} & -- & -- & $-784$ [$-966$, $-599$]\\
LEAR-type LASSO & 6.06 & 2{,}217 & 2{,}321 & -- & $+63$ [$-131$, $+259$] & $-721$ [$-1016$, $-434$]\\
Gradient boosting & 8.03 & 2{,}456 & 2{,}731 & -- & $+474$ [$+305$, $+638$] & $-310$ [$-548$, $-86$]\\
Regression & 8.22 & 2{,}631 & 3{,}041 & -- & $+784$ [$+599$, $+966$] & --\\
Seasonal naive & 13.17 & 9{,}394 & 7{,}849 & -- & $+5592$ [$+4664$, $+6568$] & $+4808$ [$+3870$, $+5769$]\\\midrule
V1 SPO+, unconstrained & 8.72 & 2{,}151 & 2{,}689 & 2{,}657--2{,}726 & $+432$ [$+241$, $+627$] & $-352$ [$-479$, $-229$]\\
V3 SPO+, anchored & 8.24 & 2{,}155 & 2{,}702 & 2{,}686--2{,}734 & $+445$ [$+240$, $+655$] & $-339$ [$-468$, $-214$]\\
V4 MLP SPO+, anchored & 8.10 & 2{,}203 & 2{,}425 & 2{,}320--2{,}523 & $+167$ [$-9$, $+346$] & $-617$ [$-796$, $-428$]\\
\bottomrule\end{tabular}

\begin{tablenotes}\footnotesize
\item Oracle: EUR\,\DEoracleH\ per MW-year over the holdout. The price limits never bind, so V2 is the
same as V1.
\end{tablenotes}
\end{threeparttable}
\end{table}

\section{Other results}\label{app:extra}
\textbf{Conditions under which averaging helps.} For any $z\in\F$, the set
$N(z)=\{\pi:z\in\argmax_{z'\in\F}J(\pi,z')\}$ equals
$\bigcap_{z'\in\mathcal V}\{\pi:\Delta\pi^\top(q(z)-q(z'))\ge C(z)-C(z')\}$, the normal cone of $z$ used in Proposition~\ref{prop:margin}(i). It is an
intersection of half-spaces, so it is convex. So, if $z$ is a best schedule for each forecast, it is also a best schedule for any weighted average of them.

\textbf{Gain from averaging and disagreement between forecasts (India, full sample).} For each day, we
measure how much the three forecasts in EW3 disagree by the average of
$\Delta\sum_t|\hat q_{k,t}-\hat q_{l,t}|$ over the three pairs. This uses forecasts only. The gain of EW3 over the mean regret of its three component forecasts, in Rs lakh per MW-year with 95\% bootstrap intervals, grows
with disagreement:
\begin{itemize}[leftmargin=1.2em]
\item low disagreement (bottom third of days): 0.23 [0.17, 0.29];
\item middle third: 0.52 [0.41, 0.63];
\item high disagreement (top third): 1.18 [0.99, 1.41].
\end{itemize}
The gain does not change in a clear way with the oracle's local dual margin (0.57, 0.71, 0.66). EW3 has lower regret than the mean regret of its three component forecasts on 75\% of days.

\textbf{Other methods tested.} Two other corrections did not help. Adjusting the trading costs online,
which by Lemma~\ref{lem:cons} is the same as rescaling and shifting the forecast, lowers regret by at most
0.07 lakh, even with hindsight. Choosing the averaging weights online from past profit or past squared
error does not outperform equal weights (full sample; differences $+0.04$ [$-0.01$, $+0.09$] and $+0.01$
[$-0.03$, $+0.04$]).

\paragraph{The 2023 cap change.} The cap fell from Rs\,12{,}000 to Rs\,10{,}000/MWh on 4~April~2023, which
changes how often prices reach the cap. Over the 90 days before and after the change, the share of realized
blocks at the cap rose from 13.1\% to 18.1\% (8.6\% and 17.2\% over all days under each cap). The cost of
projection did not follow the cap. For V1 (seeds 1 and 2), the share of forecast blocks above the cap fell
from 15.0\% to 3.1\% and the cost of projection from $+3.02$ to $-0.05$ lakh per MW-year, because its early
level drift faded as training data accumulated. For the statistical forecasts and the anchored V3, clipping
rose somewhat after the change (to at most 9.1\% of blocks), but the cost of projection stayed within $\pm0.1$
lakh per MW-year in both windows. Because the change almost coincides with the end of the first year, these
comparisons are descriptive; they do not identify an effect of the cap. They are, however, consistent with
the theory: projection is costly when the forecast level is too high, not just because the cap binds more often.

\section{Operating rules: details}\label{app:rules}
Section~\ref{sec:rule} summarizes these comparisons; this section gives the full results. Table~\ref{tab:policy}
compares always clipping (P), never clipping (U), always anchoring before clipping (A) and the
certificate-based switch (C), which clips when the conditions of Theorem~\ref{thm:safe} hold and otherwise
anchors before clipping. Table~\ref{tab:select} compares the solver's own choice among tied schedules (P)
with schedule-preserving selection (S, Corollary~\ref{cor:select}) and lexicographic selection (L), which on
the days when the raw schedule is no longer optimal picks, among the schedules that are optimal at the clipped
forecast, one that maximizes profit at the raw forecast. Rules C, S and L were each fixed before their
comparison, and remaining ties under L are broken in favor of the raw schedule.

\begin{table}[!htbp]\centering\small
\begin{threeparttable}
\caption{Ways of using a clipped forecast: regret (Rs lakh per MW-year). P: always clip; U: never clip;
A: always anchor to the regression's level, then clip; C: clip if the certificate of Theorem~\ref{thm:safe}
holds, otherwise anchor, then clip.}\label{tab:policy}
\footnotesize\setlength{\tabcolsep}{3pt}
\begin{tabular}{lccccccccc}\toprule
& & \multicolumn{4}{c}{Validation year} & \multicolumn{2}{c}{Holdout} & \multicolumn{2}{c}{Fresh test}\\
\cmidrule(lr){3-6}\cmidrule(lr){7-8}\cmidrule(lr){9-10}
Method & \makecell{Certificate\\holds (\%)} & P & U & A & C & P & C $-$ P & P & C $-$ P [95\% CI]\\\midrule
Regression & 100 & 4.45 & 4.45 & 4.45 & 4.45 & 3.17 & $0.00$ & 2.67 & $0.00$ [$0.00$, $0.00$]\\
LEAR-type LASSO & 92 & 5.24 & 5.24 & 5.24 & 5.24 & 2.71 & $0.00$ & 2.90 & $+0.67$ [$+0.04$, $+1.53$]\\
Gradient boosting & 94 & 5.69 & 5.66 & 5.74 & 5.69 & 2.94 & $0.00$ & 3.70 & $0.00$ [$0.00$, $0.00$]\\
EW3 & 100 & 4.08 & 4.08 & 4.05 & 4.08 & 2.43 & $0.00$ & 2.90 & $0.00$ [$0.00$, $0.00$]\\\midrule
V1 SPO+, unconstrained & 74 & 6.55 & 5.68 & 4.95 & 5.62 & 3.09 & $+0.02$ & 3.13 & $-0.13$ [$-0.32$, $+0.02$]\\
V2 SPO+, after projection & 65 & 6.95 & 5.80 & 5.28 & 5.72 & 3.14 & $-0.01$ & 3.46 & $-0.50$ [$-1.27$, $+0.16$]\\
V3 SPO+, anchored & 75 & 4.86 & 4.87 & 4.83 & 4.86 & 3.10 & $0.00$ & 3.12 & $0.00$ [$0.00$, $0.00$]\\
V4 neural-network SPO+ & 83 & 6.61 & 6.58 & 6.61 & 6.60 & 2.89 & $0.00$ & 3.32 & $0.00$ [$0.00$, $0.00$]\\
V5 Fenchel--Young & 90 & 5.56 & 5.56 & 5.53 & 5.56 & 3.50 & $0.00$ & 4.65 & $0.00$ [$0.00$, $0.00$]\\
\bottomrule\end{tabular}

\begin{tablenotes}\footnotesize
\item The rule was fixed before evaluation. DFL rows average three seeds. Tables~\ref{tab:policy}
and~\ref{tab:select} use the forecasts of the rerun on the extended data (Section~\ref{sec:fresh}). These
equal the forecasts behind Table~\ref{tab:iex} for every method except the neural network V4, whose rerun
differs slightly; its holdout regret is 2.894 here and 2.895 in Table~\ref{tab:iex}, which round to 2.89 and
2.90. ``Certificate holds'': share of all
evaluation days on which the conditions of Theorem~\ref{thm:safe} hold. Intervals: circular 7-day block
bootstrap. On the holdout, all differences C $-$ P lie between $-0.01$ and $+0.02$, with intervals within
$\pm0.05$.
\end{tablenotes}
\end{threeparttable}
\end{table}

\begin{table}[!htbp]\centering\small
\begin{threeparttable}
\caption{Selection among tied schedules after clipping: the solver's own choice (P), schedule-preserving
selection (S) and lexicographic selection (L). Regret in Rs lakh per MW-year.}\label{tab:select}
\footnotesize\setlength{\tabcolsep}{2.4pt}
\begin{tabular}{lcccccccc}\toprule
& & \multicolumn{3}{c}{Validation year} & \multicolumn{2}{c}{Holdout} & \multicolumn{2}{c}{Fresh test}\\
\cmidrule(lr){3-5}\cmidrule(lr){6-7}\cmidrule(lr){8-9}
Method & \makecell{Raw schedule optimal\\after projection (\%)} & P & S & L & P & S $-$ P [95\% CI] & P & S $-$ P [95\% CI]\\\midrule
Regression & 100.0 & 4.45 & 4.45 & 4.45 & 3.17 & $0.00$ [$0.00$, $0.00$] & 2.67 & $0.00$ [$0.00$, $0.00$]\\
LEAR-type LASSO & 99.3 & 5.24 & 5.24 & 5.24 & 2.71 & $0.00$ [$-0.01$, $0.00$] & 2.90 & $0.00$ [$0.00$, $0.00$]\\
Gradient boosting & 100.0 & 5.69 & 5.66 & 5.66 & 2.94 & $0.00$ [$-0.01$, $+0.01$] & 3.70 & $0.00$ [$0.00$, $0.00$]\\
EW3 & 100.0 & 4.08 & 4.08 & 4.08 & 2.43 & $0.00$ [$0.00$, $0.00$] & 2.90 & $0.00$ [$0.00$, $0.00$]\\\midrule
V1 SPO+, unconstrained & 98.3 & 6.55 & 6.50 & 6.15 & 3.09 & $+0.02$ [$0.00$, $+0.05$] & 3.13 & $-0.15$ [$-0.33$, $0.00$]\\
V2 SPO+, after projection & 97.8 & 6.95 & 6.82 & 6.40 & 3.14 & $+0.03$ [$-0.02$, $+0.09$] & 3.46 & $-0.52$ [$-1.25$, $+0.09$]\\
V3 SPO+, anchored & 100.0 & 4.86 & 4.87 & 4.87 & 3.10 & $+0.03$ [$+0.01$, $+0.07$] & 3.12 & $0.00$ [$0.00$, $0.00$]\\
V4 neural-network SPO+ & 99.9 & 6.61 & 6.58 & 6.58 & 2.89 & $-0.04$ [$-0.06$, $-0.02$] & 3.32 & $0.00$ [$0.00$, $0.00$]\\
V5 Fenchel--Young & 99.6 & 5.56 & 5.56 & 5.56 & 3.50 & $0.00$ [$0.00$, $0.00$] & 4.65 & $0.00$ [$0.00$, $0.00$]\\
\bottomrule\end{tabular}

\begin{tablenotes}\footnotesize
\item Rules S and L fixed before evaluation. As in Table~\ref{tab:policy}, the V4 row uses the rerun
of the neural network (holdout regret 2.894, against 2.895 in Table~\ref{tab:iex}). L differs from S only on days when the raw schedule is not optimal
after clipping; on the holdout and fresh days, L and S give the same regret to two decimals for every method.
DFL rows average three seeds. The second column is the share of all
evaluation days on which the raw schedule is a best schedule at the clipped forecast. Intervals: circular
7-day block bootstrap.
\end{tablenotes}
\end{threeparttable}
\end{table}

\end{document}